\documentclass{article}
\usepackage{arxiv}
\usepackage{float}
\usepackage{algorithm,algpseudocode}
\usepackage{amsmath,amssymb}
\usepackage{amsthm,latexsym,float}
\usepackage[labelformat=simple]{subfig}

\usepackage{bm,multirow}        
\usepackage{graphicx}
\usepackage{xcolor}
\usepackage{placeins}
\usepackage{doi}

\newcommand{\bh}{\mathbf{h}}
\newcommand{\bu}{\mathbf{u}}
\newcommand{\bv}{\mathbf{v}}
\newcommand{\bs}{\mathbf{s}}
\newcommand{\bb}{\mathbf{b}}
\newcommand{\bw}{\mathbf{w}}
\newcommand{\bF}{\mathbf{F}}
\newcommand{\bA}{\mathbf{A}}
\newcommand{\bH}{\mathbf{H}}
\newcommand{\bI}{\mathbf{I}}
\newcommand{\cD}{\mathcal{D}}
\newcommand{\cG}{\mathcal{G}}
\newcommand{\cX}{\mathcal{X}}
\newcommand{\cA}{\mathcal{A}}
\newcommand{\cH}{\mathcal{H}}
\newcommand{\cL}{\mathcal{L}}

\usepackage{hyperref}

\title{Intrusive and non-intrusive reduced order modeling of the rotating thermal shallow water equation}

\author{ \href{https://orcid.org/0000-0001-7904-605X}{\includegraphics[scale=0.06]{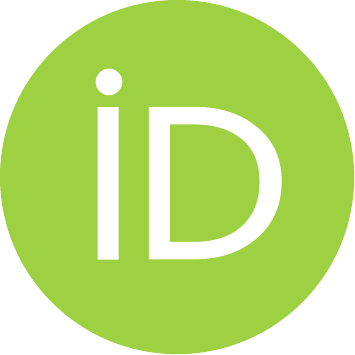}\hspace{1mm}S\"uleyman Y{\i}ld{\i}z} \\
Institute of Applied Mathematics\\
Middle East Technical University, 
 Ankara-Turkey\\
	\texttt{yildiz.suleyman@metu.edu.tr} \\
	\And
	\href{https://orcid.org/0000-0003-1037-5431}{\includegraphics[scale=0.06]{orcid.pdf}\hspace{1mm}B\"ulent Karas\"ozen} \\
     Institute of Applied Mathematics \& Department of Mathematics\\
     Middle East Technical University,
     Ankara-Turkey\\
     \texttt{bulent@metu.edu.tr}
	\And
  \href{https://orcid.org/0000-0001-5262-063X}{\includegraphics[scale=0.06]{orcid.pdf}\hspace{1mm}Murat Uzunca} \\
   Department of Mathematics\\
	Sinop University,
     Sinop-Turkey \\
     \texttt{muzunca@sinop.edu.tr}\\
}

\date{}

\begin{document}
\maketitle

\begin{abstract}
In this paper, we investigate projection-based intrusive and data-driven non-intrusive model order reduction methods in numerical simulation of rotating thermal shallow water equation (RTSWE) in parametric and non-parametric form. Discretization of the RTSWE in space with centered finite differences leads to Hamiltonian system of ordinary differential equations with linear and quadratic terms. The full-order model (FOM) is obtained by applying linearly implicit Kahan's method in time. Applying proper orthogonal decomposition with Galerkin projection (POD-G), we construct the intrusive reduced-order model (ROM). We apply operator inference (OpInf) with re-projection for non-intrusive reduced-order modeling. In the parametric case, we make use of the parameter dependency at the level of the PDE without interpolating between the reduced operators. The least-squares problem of the OpInf is regularized with the minimum norm solution. Both ROMs behave similar and are able to accurately predict the test and training data and capture system behavior in the prediction phase with several orders of computational speedup over the FOM. The preservation of system physics such as the conserved quantities of the RTSWE by both ROMs enables that the models fit better to data and stable solutions are obtained in long-term predictions, which are robust to parameter changes.
\end{abstract}


\section{Introduction}

Rotating shallow water equation (RSWE) \cite{Salmon04} is used as  conceptual and predictive model in geophysical fluid dynamics for the behavior of rotating inviscid fluids and  has  been applied to
a diverse range of oceanographic and atmospheric phenomena. To
describe three dimensional flows, 
multilayer RSWE model is used, where the fluid is discretized horizontally into a
set of vertical layers, each having its own height, density and horizontal velocity.  
However, the RSWE model does not allow for gradients of the mean temperature and/or density. The rotating thermal shallow water equation (RTSWE)  \cite{Dellar03,Dellar13,Eldred19,Ripa95}  represents an extension of the RSWE equation, to include horizontal density/temperature gradients both in the atmospheric and oceanic context. The  RTSWE is used in the general circulation models \cite{Zerroukat15}, planetary flows \cite{Dellar14},  and modeling atmospheric and oceanic temperature fronts \cite{Dempsey88,Young95}, thermal instabilities \cite{Gouzien17}.

High fidelity full order models (FOMs) for  partial differential equations (PDEs) such as the shallow water equations (SWEs)  are obtained  with the standard discretization methods, e.g., finite-differences, finite-volumes, finite-elements, spectral elements, or discontinuous Galerkin methods. The computational cost associated with fully resolved solutions remains a barrier in
real-time simulations on fine space-time grids.
The dynamics of large-scale dynamical systems such as the SWEs typically lie in a low-dimensional space.
Model  order  reduction (MOR)  aims to construct low-dimensional and efficient reduced-order models (ROMs) to accurately approximate the underlying FOMs, which requires repeated evaluations of the model over a range of parameter values. Additionally, ROMs are even more valuable for SWEs in simulating and predicting the model over a  long time horizon.The ROMs are implemented using the offline-online approach.  The reduced basis functions are extracted from the snapshots in the offline stage,  and the reduced basis is computed by combining them. In the online phase, the full problem is projected onto the reduced space, and the solutions for new parameters are computed in an efficient manner. Based on the offline-online methodology, there exists two versions of ROM methods; the intrusive and non-intrusive ones.
In the case of intrusive ROM methods, the reduced solutions are determined by solving a reduced order model, i.e., a projection of the FOM onto the reduced space.
The proper orthogonal decomposition (POD) with the Galerkin projection (POD-G) \cite{Berkooz93,Sirovich87} is arguably one of the most popular intrusive ROM methods.
Applying the singular value decomposition (SVD) to the snapshot matrix, the POD basis are extracted.  Then, ROM is constructed by applying Galerkin projection. For PDEs with only linear and polynomial terms like the SWEs,  the projection-based reduced model admits an efficient low-dimensional numerical representation independent of the dimension of the FOM without necessitating hyper-reduction techniques like empirical/discrete empirical interpolation method \cite{Barrault04,chaturantabut10nmr}. We refer to the books \cite{Benner17,Rozza14,Benner21a} for an overview of the available MOR techniques.

Traditionally, projection-based model reduction is intrusive since the numerical implementation of the reduced models require access to the discretized PDE operators.
The intrusive nature of MOR techniques limits the scope of traditional model reduction methods. The major drawback of intrusive methods is that they require access to the full model and solvers are typically unavailable when working with, and thus the traditional intrusive model reduction methods are not applicable when  proprietary software is used for solving the PDEs, where the details of the governing equations are not known. 
Another class of ROMs is the data-driven or non-intrusive ROMs, which are fundamentally different from the intrusive ROMs. Unlike intrusive, using the non-intrusive model reduction techniques, reduced models are learnt from snapshots, i.e., either numerical approximations or measurements of the states of the dynamical systems, when the operators of the discretized systems are not available. There exist several software packages that are able to simulate SWEs for a given parameter set and an initial condition \cite{Cotter19b,Delestre17}.

Machine learning plays an important role in analyzing the underlying process of the dynamics from the data. Recent advances in machine learning techniques such as neural networks offer new opportunities to develop more efficient and accurate ROMs. They learn a model based on the training data that neither requires explicit access to the high-fidelity model operators nor any additional information about the process. However, the amount of data required to learn the model accurately imposes a burden in the context of large-scale PDE simulations \cite{swischuk2019projection}. Incorporating existing knowledge about the physics of the models requires less training data.
We review the most relevant literature  about learning dynamical-system models from data  related to this work.
Data-driven  reduced models are constructed in the Loewner framework from input-output measurements by fitting linear-time invariant models to frequency or time-domain data \cite{Antoulas14,Peherstorfer17}.
The Loewner approach has been extended to bilinear \cite{Antoulas16} and quadratic-bilinear systems \cite{Antoulas18}.
Dynamic mode decomposition (DMD) learns reduced models for nonlinear dynamical systems, by fitting linear operators to state trajectories with respect to the $L_2$ norm  \cite{Rowley09sf,Schmid10dmd}. Methods based on the Koopman operator has been developed to extend DMD to nonlinear dynamical systems \cite{Williams15}.
Reduced models are also constructed  in the high-dimensional systems by exploiting sparsity  \cite{Brunton16a,Loiseau18,Schaeffer18}.
Recently, data-fit surrogate models with artificial neural networks (ANNs) are used as the regression models for the non-intrusive ROMs for time-dependent dynamical systems \cite{HesthavenUbbiali18,WangHesthaven19}.

In recent years, the operator inference (OpInf) method to construct ROMs non-intrusively has gained much attention, which was first investigated for polynomial nonlinearities in \cite{Peherstorfer16}. The operators defining the ROM can be learned by formulating an optimization problem, i.e., least-squares problem without accessing the discretized operators of the PDEs.
The methodology was later extended to nonlinear systems that can be written as a polynomial or quadratic-bilinear system by utilizing lifting variable transformation, called Lift \& Learn method  \cite{Kramer19,Kramer20}. In \cite{Kramer20a}, OpInf method is applied to nonlinear systems in which the structure of the nonlinearities is preserved while learning ROMs from data. Data-driven ROMs via regularized  OpInf are employed to combustion problems \cite{Willcox20,Swischuk20}. The OpInf method is also used for data-driven reduced-order modeling in fluid dynamics; such as incompressible flows \cite{Benner20} and SWEs \cite{Karasozen20a}. Recently in \cite{Willcox21}, the OpInf method is generalized to the PDE setting through the use of lifting variables. In the case of linear FOMs, the  OpInf method equivalent to the DMD, the inferred operators are the same. 
Recently the DMD is extended to quadratic bilinear systems \cite{gosea2020fitting}. When the future states of dynamical systems depend only on the current state and not on previous ones, they are called Markovain systems. The non-Markovian systems are considered having a memory, where future states depend on the current and previous states.
The projected trajectories   of   the OpInf \cite{Peherstorfer16}  correspond to non-Markovian dynamics in the low-dimensional subspaces even though the FOMs are Markovian, as  known from, e.g., the Mori--Zwanzig formalism \cite{Chorin06,Stuart04}.
The non-Markovian dynamics has been investigated intensively \cite{Duraisamy17,Iliescu21,Maday18,San19a,Wang20}. To recover the Markovian dynamics of the reduced system, a data sampling scheme has devised in \cite{Peherstorfer20}, which iterates between time stepping the high-dimensional FOM and projecting onto low-dimensional reduced subspaces.  In this way the  Markovian dynamics is retained in low-dimensional subspaces. It was shown in \cite{Peherstorfer20} that under certain conditions, applying operator inference to re-projected trajectories for dynamical systems with polynomial nonlinear terms, produces the same operators that are obtained with the intrusive model reduction methods.
Recently, probabilistic a posteriori error estimators are derived for the OpInf with re-projection for linear parabolic PDEs \cite{Peherstorfer20a} and a non-Markovian OpInf with partial information is investigated  \cite{Peherstorfer21}. A deep learning version of the OpInf is introduced in \cite{Benner21deep}.

Intrusive and non-intrusive MOR techniques, POD-G and DMD  have been intensively studied for the SWEs, see e.g., \cite{Bistrian15,Bistrian17,Esfahanian09,Karasozen21,Lozovskiy16,Lozovskiy17,Navon20,cstefuanescu2014comparison}.
In this paper, we present efficient ROMs by applying the intrusive POD-G and the non-intrusive data-driven OpInf method to the RTSWE in the parametric and non-parametric forms. The RTSWE is discretized in space with the second-order centered finite differences while preserving the  Hamiltonian structure, which leads to a linear-quadratic system of ordinary differential equations (ODEs). For time discretization we use the linearly implicit Kahan's method which is designed for linear-quadratic ODE systems like the RTSWE in the semi-discrete form. It is second-order accurate in time, and has a lower computational cost than the fully implicit schemes such as the Crank-Nicolson scheme. Kahan's method requires only one step Newton iteration at each time level.  Furthermore, it has favorable conservation properties; the energy (Hamiltonian) and other conserved quantities of the RTSWE like the mass, the total vorticity, and the buoyancy are preserved in long-time integration.
To accelerate the online and offline ROM computations, we apply the POD-G  by the use of matricization of tensors  \cite{Benner18,Benner21,Karasozen21} and the sparse matrix technique MULTIPROD \cite{leva08mmm}.
We apply  the OpInf with re-projection \cite{Peherstorfer20} to the RTSWE by preserving the low-dimensional Markovian dynamics.  It was shown through by the numerical results that learning low-dimensional models with the OpInf method with re-projection \cite{Peherstorfer20} are very close to the reduced models from intrusive model reduction POD-G. The energy and other conserved quantities are well preserved with the POD-G and OpInf without any drift over time.
We compare POD-G with the  OpInf method for the parametric case with the Coriolis force, where we make use of the known parametric dependency at the PDE level, whereas in \cite{Peherstorfer16,Peherstorfer20} the ROMs are constructed at each training parameter via interpolation. We show that under certain assumptions on the time discretization \cite{Peherstorfer16,Peherstorfer20} the intrusive POD-G model converges to the learned non-intrusive OpInf model.
Numerical results also show that both the POD-G and the OpInf methods can predict the reduced parametric models with high accuracy. Furthermore, we show that the ROMs exhibit high accuracy in the training regime and acceptable accuracy in prediction phase. The data matrices in the least-squares problem of the OpInf have large condition numbers, leading to ill-conditioned inverse problem. To deal with this issue, regularization techniques  like the truncated SVD or Tikhanov regularization should be used. Because, decay of the singular values of the data matrices does not provide any information about the choice of the tolerance for the regularization parameter, we solve the least-squares problem in the minimum norm, where the tolerances are determined with the L-curve. Speed-up factors of order two over the FOM are achieved for both ROMs, whereas the  OpInf is more costly than the POD-G due to the re-projection.

The paper is organized as follows. In Section~\ref{sec:fom}, the RTSWE and the full space-time discretization are described. The detailed formulations of intrusive POD-G and the non-intrusive  OpInf are given in Section~\ref{sec:rom}. In Section~\ref{sec:num}, numerical results are presented for the double vortex RTSWE in the non-parametric and in the parametric form, and for prediction over a time horizon. The paper ends with some concluding remarks in Section~\ref{sec:conc}.

\section{Rotating Thermal Shallow Water Equation}
\label{sec:fom}

RSWE \cite{Salmon04,Salmon07} is considered in rotating inviscid fluids with one or more layers \cite{Stewart16}, within each layer the horizontal velocity is assumed to be depth-independent, so the fluid moves in columns. Large-scale vortex and wave dynamics can be understood with this model \cite{Zeitlin18}. However, a drawback of this model is that it does not admit horizontal temperature and density gradients. The RTSWE represents an extension of the RSWE  allowing for  gradients  in the atmospheric and oceanic context, known also as Ripa equation \cite{Ripa95}. It is obtained along the same lines as the RSWE, using the hypothesis of columnar motion, but with nonuniform density/temperature \cite{Gouzien17,Zeitlin18}.
Typical applications of the RTSWE  are to model atmospheric and oceanic temperature fronts, mixed boundary layers \cite{Young95}.
RTSWE is given for the primitive variables \cite{Dellar13,Eldred19} as
\begin{equation} \label{tswe}
\begin{aligned}
\frac{\partial h}{\partial t } & =   - (hu)_x   - (hv)_y,\\
\frac{\partial u}{\partial t } & =  hqv -\left (\frac{u^2+v^2}{2}\right)_x- \frac{h}{2}s_x - s( h + b)_x,   \\
\frac{\partial v}{\partial t } & =   -hqu -\left (\frac{u^2+v^2}{2}\right)_y - \frac{h}{2}s_y - s(h + b  )_y, \\
\frac{\partial s}{\partial t } & =   - us_x   - vs_y,
\end{aligned}
\end{equation}
where $u(x,y,t)$ and $v(x,y,t)$ are the relative velocities, $h(x,y,t)$ is the fluid height, $b(x,y)$ topographic height, $\rho(x,y,t)$ fluid density, $g$ gravity constant, $s = g\frac{\rho}{\bar{rho}}$ the buoyancy, $S= hs$ mass-weighted buoyancy, and
$q = (v_x - u_y + f)/h$ is the potential vorticity with the Coriolis parameter $f(\mu) = 2\Omega\sin(\mu)$,  where $\mu $ varies with the latitude and $\Omega = 7.292\times 10^{-5}$ is the rotational speed of the Earth.
The standard RSWE is recovered in the limit of constant density, i.e., constant buoyancy.

The RTSWE \eqref{tswe} is considered on a rectangular domain $\Omega = [a,b] \times [c,d]$ with  periodic boundary conditions, and on a time interval $[0,T]$ for a target time $T>0$, together with the given initial conditions
$$
u({\bm x},0) = u_0({\bm x}), \; v({\bm x},0) = v_0({\bm x}), \; h({\bm x},0) = h_0({\bm x}), \;
s({\bm x},0) = s_0({\bm x}).
$$

The RTSWE is a non-canonical Hamiltonian PDE with state-dependent Poisson matrix \cite{Dellar02,Dellar13,Eldred19}. The conserved Hamiltonian or energy is given by
\begin{equation} \label{ham}
\mathcal{E}(z)= \int_{\Omega}  \left (\frac{h^2s}{2} + hsb + h \frac{u^2 + v^2}{2} \right ) d\Omega.
\end{equation}

Other conserved quantities are the Casimirs: the  mass, the total potential vorticity, and  the buoyancy defined respectively by
\begin{equation} \label{cas}
{\mathcal M} =\int_{\Omega} h \;d\Omega, \quad {\mathcal Q }=\int_{\Omega} q \;d\Omega, \quad {\mathcal B} =\int_{\Omega} hs \; d\Omega.
\end{equation}

The RTSWE \eqref{tswe} is discretized by finite differences on a uniform grid of the rectangular spatial domain $\Omega =[a,b]\times [c,d]$ with the nodes ${\bm x}_{ij} = (x_i,y_j )^T$, where $x_i=a+(i-1)\Delta x$, $y_j=c+(j-1)\Delta y$, $i,j=1,\ldots, n+1$, and spatial mesh sizes are $\Delta x=(b-a)/n$ and $\Delta y=(d-c)/n$. Using the from bottom to top/from left to right ordering of the grid nodes, the vectors of semi-discrete state variables are defined as
\begin{equation}\label{solvec}
\begin{aligned}
\bm{h}(t) &= (h_{1,1}(t),\ldots , h_{1,n}(t),h_{2,1}(t),\ldots , h_{2,n}(t), \ldots , h_{n,n}(t))^T, \\
\bm{u}(t) &= (u_{1,1}(t),\ldots , u_{1,n}(t),u_{2,1}(t),\ldots , u_{2,n}(t), \ldots , u_{n,n}(t))^T,\\
{\bm v}(t) &= (v_{1,1}(t),\ldots , v_{1,n}(t),v_{2,1}(t),\ldots , v_{2,n}(t), \ldots , v_{n,n}(t))^T,\\
{\bm s}(t) &= (s_{1,1}(t),\ldots , s_{1,n}(t),s_{2,1}(t),\ldots , s_{2,n}(t), \ldots , s_{n,n}(t))^T,
\end{aligned}
\end{equation}
where we set $w_{i,j}(t)=w(x_i,y_j,t)$ for $w=h,u,v,s$.
We note that the degree of freedom (dimension of FOM) is given by $N:=n^2$ because of the periodic boundary conditions, i.e., the most right grid nodes ${\bm x}_{n+1,j}$ and the most top grid nodes ${\bm x}_{i,n+1}$ are not included. The solution vector is defined by ${\bm z}(t)=(\bm{h}(t),\bm{u}(t),{\bm v}(t), {\bm s}(t)): [0,T]\mapsto \mathbb{R}^{4N}$.
Throughout the paper, we do not explicitly represent the time dependency of the semi-discrete solutions for simplicity, and we write $\bm{u}$, ${\bm v}$, $\bm{h}$, ${\bm s}$, and  ${\bm z}$.

The first-order derivatives in space are approximated by utilizing one-dimensional centered finite differences in $x$ and $y$ directions, and they are extended to two dimensions by the use of the Kronecker product. Let $\widetilde{D}_n\in \mathbb{R}^{n\times n}$ be the matrix corresponding to the centered finite differences under periodic boundary conditions on a one-dimensional domain (interval)
$$
 \widetilde{D}_n=
\begin{pmatrix}
 0& 1&  & &-1 \\
-1& 0&1 & &   \\
  & \ddots  & \ddots  &\ddots  &   \\
  &  &-1&0 &1 \\
 1&  &  & 1&0
\end{pmatrix}.
$$
Then, on the two dimensional mesh, the centered finite difference matrices $D_x\in\mathbb{R}^{N\times N}$ and $D_y\in\mathbb{R}^{N\times N}$ corresponding to the first order partial derivatives $\partial_x$ and $\partial_y$, respectively, can be calculated as
$$D_x=\frac{1}{2\Delta x}\widetilde{D}_{n}\otimes I_{n} \; ,  \quad D_y=\frac{1}{2\Delta y}I_{n}\otimes \widetilde{D}_{n},
$$
where $\otimes$ denotes the Kronecker product, and $I_n$ is the identity matrix of size $n$.
Then, the semi-discretization of the RTSWE  \eqref{tswe} leads to the following system of linear-quadratic ODEs
\begin{equation}\label{swef}
\begin{aligned}
\frac{d\bh}{dt }  &=   -D_x(\bu\circ \bh)-D_y(\bv\circ \bh),\\
\frac{d\bu}{dt }  &= - \bu\circ (D_x\bu) - \bv\circ( D_y\bu)
 - \frac{\bh}{2} \circ  D_x\bs   -  \bs \circ  D_x\bh -
  \bs\circ  D_x\bb    +   f \bv, \\
\frac{d\bv}{dt }  &=  -\bu\circ( D_x\bv)  - \bv\circ( D_y\bv)
 - \frac{\bh}{2} \circ  D_y\bs   -  \bs \circ  D_y\bh -  \bs\circ  D_y\bb
- f \bu, \\
\frac{d\bs}{dt } & = - \bu \circ D_x\bs -   \bv \circ D_y\bs.
\end{aligned}
\end{equation}
For $\bw = (\bh,\bu,\bv,\bs)$, the system \eqref{swef} can be rewritten in the compact form as
\begin{align} \label{tenfom}
\dot\bw = \bF(\bw)=  \bA\bw+  \bH(\bw\otimes \bw),
\end{align}
where  $\dot{\bw} $ denotes the time derivative of the state vector $\bw$, and $\bA\in \mathbb{R}^{4N\times 4N}$ is a linear matrix operator corresponding to the linear terms
\begin{equation*}
\quad \bA =
\begin{pmatrix}
0 & 0 & 0 & 0 \\
 0 & 0 & f\bI & - (D_x\bb)^d \\
 0 &- f \bI & 0 & - (D_y\bb)^d \\
0 & 0 & 0 & 0
\end{pmatrix},
\end{equation*}
where $ \bI \in \mathbb{R}^{N\times N} $ is the identity matrix and $(D_x\bb)^d\in \mathbb{R}^{N\times N} $ is the diagonal matrix whose diagonal entries are given by the elements of the vector $D_x\bb$, i.e., $((D_x\bb)^d)_{ii}=(D_x\bb)_i$, $i=1,\ldots ,N$.
The vector $\bH(\bw\otimes \bw)\in \mathbb{R}^{4N}$ contains the quadratic terms
\begin{equation*}
\bH(\bw\otimes \bw)=
\begin{pmatrix}
- D_x (\bu\circ \bh)  - D_y (\bv\circ \bh)\\
- \bu\circ (D_x\bu) - \bv\circ (D_y\bu) - \frac{\bh}{2} \circ  D_x\bs   -  \bs \circ  D_x\bh  \\
- \bu\circ (D_x\bv)  -  \bv\circ (D_y\bv)  - \frac{\bh}{2} \circ  D_y\bs   -  \bs \circ  D_y\bh  \\
- \bu \circ D_x\bs -   \bv \circ D_y\bs
\end{pmatrix},
\end{equation*}
where $\circ$ denotes the elementwise matrix or Hadamard product, and $\bH\in \mathbb{R}^{4N\times (4N)^2}$ is an appropriate matricized tensor operator such that the right hand side of the above equality is obtained.

For time discretization, we partition the time interval $[0,T]$ into $K$ uniform intervals with the step size $\Delta t=T/K$ as $0=t_0<t_1<\ldots <t_{K}=T$, and $t_k=k\Delta t$, $k=0,1,\ldots ,K$. Then, we denote by $\bw^k=\bw(t_k)$ the fully discrete solution vector at time $t_k$. For the linear-quadratic autonomous ODE system \eqref{tenfom},  we use Kahan's "unconventional"  discretization \cite{Kahan93}, which yields the system
\begin{equation*}
\frac{\bw^{k+1} - \bw^k}{\Delta t} = \frac{1}{2}\bA\left(\bw^k + \bw^{k+1}\right) +
 \overline{\bH}\left((\bw^k\otimes \bw^k),(\bw^{k+1}\otimes \bw^{k+1})\right),
\end{equation*}
where the symmetric bilinear form $\overline{\bH}(\cdot ,\cdot)$  denotes  the polarization \cite{Celledoni15}  of the quadratic vector field $\bH(\cdot )$, given by
\begin{align*}
\overline{\bH}\left((\bw^k\otimes \bw^k),(\bw^{k+1}\otimes \bw^{k+1})\right) =& \;
\frac{1}{2}\left(\bH\left((\bw^k+\bw^{k+1})\otimes(\bw^k+\bw^{k+1})\right)\right)   \\ & -  \frac{1}{2}\left(\bH\left(\bw^k\otimes \bw^k\right) +
 \bH\left(\bw^{k+1}\otimes \bw^{k+1}\right)\right).
\end{align*}
Kahan's method is second-order and time-reversal and linearly implicit for quadratic vector fields  \cite{Celledoni13}, i.e., it requires only one step of Newton iterations per time step
\begin{equation} \label{kahan}
\left ( I -\frac{\Delta t}{2} \bF'(\bw^k)\right ) \frac{ \bw^{k+1} - \bw^k}{\Delta t}  =  \bF(\bw^k),
\end{equation}
where $\bF'(\bw^k)$ stands for the Jacobian matrix of $\bF(\bw)$ evaluated at $\bw^k$, and $I$ is the identity matrix of size $4N$.

Under periodic boundary conditions, the full discrete forms of the energy, the mass, the total potential vorticity, and the buoyancy are given at the time instance $t_k$ respectively as
\begin{equation}\label{dener}
\begin{aligned}
H^k(\bw) &= \sum_{i=1}^n\sum_{j=1}^n \left (\frac{1}{2}(h_{i,j}^k)^2 s_{i,j}^k  + h_{i,j}^k s_{i,j}^k b_{i,j}^k + h_{i,j}^k\frac{(u_{i,j}^k)^2 + (v_{i,j}^k)^2}{2}   \right ) \Delta x\Delta y,
\\
M^k(\bw) &= \sum_{i=1}^n\sum_{j=1}^n  h_{i,j}^k \Delta x\Delta y, \\
Q^k(\bw) & =  \sum_{i=1}^n\sum_{j=1}^n \left ( \frac{v_{i+1,j}^k-v_{i-1,j}^k}{2\Delta x} -  \frac{u_{i,j+1}^k-u_{i,j-1}^k}{2\Delta y} + f \right ) \Delta x\Delta y, \\
B^k(\bw) &= \sum_{i=1}^n\sum_{j=1}^n h_{i,j}^k s_{i,j}^k \Delta x\Delta y,
\end{aligned}
\end{equation}
where $b_{i,j}^k= b(x_i,y_j,t_k)$, and, for $w=h,u,v,s$, $w_{i,j}^k= h(x_i,y_j,t_k)$ are full discrete approximations at $(x_i,y_j,t_k)$ satisfying  $w_{i+n,j}^k=w_{i,j}^k$  and $w_{i,j+n}^k=w_{i,j}^k$  because of the periodic boundary conditions, $i,j=1,\ldots ,n$.

Linear conserved quantities such as the mass are preserved by all Runge-Kutta type methods including Kahan's methods. Kahan's methods also preserves the quadratic and cubic polynomial invariants  of  Hamiltonian PDEs approximately over long times \cite{Celledoni15,Eidnes20}. The numerical errors in the conserved quantities of the RTSWE, such as the energy, vorticity and buoyancy have bounded oscillations over time.

\section{Reduced order modeling}
\label{sec:rom}

In this section the intrusive POD-G ROM  and the non-intrusive ROM with the OpInf are constructed for the following parametric form of RTSWE \eqref{tenfom}
\begin{align} \label{tenfompar}
\dot{\bw}(t;\mu) =   \bA(\mu) {\bw}(t;\mu) +  \bH({\bw}(t;\mu)\otimes {\bw}(t;\mu)).
\end{align}
We remark that only the linear operator $\bA$ depends on the parameter $\mu\in {\mathcal D} \subset \mathbb{R}^d$, the operator $\bH$ in quadratic part is independent of the parameter.

We describe first the offline stage of constructing the POD basis, which  is common for the intrusive and non-intrusive model reduction. The first step of the offline stage represents an initialization of a training set of parameters of cardinality $\mathrm{M}$.
Then, in the second step, we query the FOM solutions for each parameter $\mu$ in the training set. To each parameter, it is associated a snapshot containing approximations of the state variables $\bm{h},\bm{u},\bm{v},\bm{s}$ at each time step. In the third step, the information provided by the snapshots is compressed through the POD.

The reduced basis vectors are obtained usually by stacking all state variables in one vector and a common reduced subspace is computed by taking the SVD of the snapshot data. Because the governing PDEs like the RTSWE \eqref{tswe} are coupled, the resulting ROMs do not preserve the coupling topological structure of the FOM \cite{Benner15c,Reiss07} and both intrusive and non-intrusive methods produce unstable ROMs \cite{Kramer20}.  In order to maintain the coupling structure in ROMs, the POD basis vectors are computed separately for each  state vector  $\bh,\bu,\bv,\bs$.
Let $ W_h(\mu),W_u(\mu),W_v(\mu),W_s(\mu)\in \mathbb{R}^{N\times K} $ be the trajectories for each state vector and parameter $\mu$
\begin{align*}
&W_h(\mu) =  \left[\bh^1(\mu) \cdots \bh^K(\mu)\right],  \quad W_u(\mu)  =  \left[\bu^1(\mu)\cdots \bu^K(\mu)\right],\\
&W_v(\mu) =  \left[\bv^1(\mu) \cdots \bv^K(\mu)\right], \quad W_s(\mu) =  \left[\bs^1(\mu) \cdots \bs^K(\mu)\right].
\end{align*}
The corresponding global snapshot matrices of the concatenated  trajectories and  parameter values $ \mu_1,\ldots, \mu_M \in \cD $ are given as
\begin{align*}
&W_{h\mu} =  \left[W_h(\mu_1) \cdots W_h(\mu_M)\right], \quad W_{u\mu} =  \left[W_u(\mu_1) \cdots W_u(\mu_M)\right],\\
&W_{v\mu} =  \left[W_v(\mu_1) \cdots W_v(\mu_M)\right], \quad W_{s\mu} =  \left[W_s(\mu_1) \cdots W_s(\mu_M)\right].
\end{align*}

The POD basis are constructed by computing the SVD of the global snapshot matrix $W_\mu \in \mathbb{R}^{4N\times MK}$ for the semidiscrete RTSW equation \eqref{tenfom} in the following form
$$
W_\mu =
\begin{pmatrix} W_{h\mu} \\ W_{u\mu} \\ W_{v\mu} \\ W_{s\mu}
\end{pmatrix}
= \begin{pmatrix} \Phi_{h} &  & &  \\  & \Phi_u & & \\  & & \Phi_v & \\  & & & \Phi_s  \end{pmatrix}
\begin{pmatrix} \Sigma_h &  & &  \\  & \Sigma_u & & \\  & & \Sigma_v & \\  & & & \Sigma_s  \end{pmatrix}
\begin{pmatrix} \Psi_h^T \\  \Psi_u^T \\ \Psi_v^T \\  \Psi_s^T  \end{pmatrix}
$$
with the orthonormal matrices $\Phi_h,\Phi_u,\Phi_v,\Phi_s \in \mathbb{R}^{N\times N}$, and $\Psi_h,\Psi_u,\Psi_v,\Psi_s \in \mathbb{R}^{MK\times N}$, and  with the diagonal matrices
$\Sigma_h,\Sigma_u,\Sigma_v,\Sigma_s \in \mathbb{R}^{N\times N}$ containing  the singular    values in descending order $\sigma_1\ge \sigma_2, \ldots$. The POD basis vectors $\Phi_{r_h}$, $\Phi_{r_u}$, $\Phi_{r_v}$ and $\Phi_{r_s}$
are given by the $r_h,r_u,r_v$ and $r_s$ leading columns of the matrices  $\Phi_h,\Phi_u,\Phi_v$ and $\Phi_s$, respectively, and yields the approximation
$$
\begin{pmatrix} \bh \\ \bu \\  \bv  \\ \bs    \end{pmatrix}
\approx \begin{pmatrix} \Phi_{r_h} &  & &  \\  & \Phi_{r_u} & & \\  & & \Phi_{r_v} & \\  & & & \Phi_{r_s}  \end{pmatrix}
\begin{pmatrix} \tilde{\bh} \\ \tilde{\bu} \\  \tilde{\bv}  \\ \tilde{\bs}    \end{pmatrix},
$$
where $\tilde{\bh},  \tilde{\bu}, \tilde{\bv}$ and $\tilde{\bs}$ are the reduced coefficients.
Although  different number of modes can be used for each state, in the present study we use equal number of modes, i.e., $r= r_u = r_v = r_h = r_s$. Note that the POD basis are  independent of the parameter $\mu \in {\mathcal D}$.

The deterministic SVD for the  snapshot matrices of size $N\times K$ has asymptotic time complexity of ${\mathcal O}(\min(NK^2,N^2K)$. We use the randomized SVD  (rSVD) \cite{Halko11,Bach19} which is faster than the deterministic SVD for large matrices and has asymptotic complexity of only $\mathcal{O}(NKr)$.

\subsection{Proper orthogonal decomposition with Galerkin projection}
\label{POD-G}

In the POD-G ROM, a reduced system is obtained by Galerkin projection of  the FOM \eqref{tenfom} onto the space spanned by the POD vectors
$$
\Phi_r =
\begin{pmatrix}
\Phi_{r_h} &   &   &   \\[0.2cm]
   & \Phi_{r_u} &   &   \\[0.2cm]
   &  & \Phi_{r_v} &   \\[0.2cm]
  &   &   & \Phi_{r_s}
\end{pmatrix}\in \mathbb{R}^{4N\times 4r}.
$$
After projection, the reduced system reads as
\begin{equation}\label{tswepfom0}
\dot{\tilde{\bw}}(t;\mu)= \widetilde{\bA}(\mu)\tilde{\bw}(t;\mu) + \widetilde{\bH}(\tilde{\bw}(t;\mu)\otimes\tilde{\bw}(t;\mu)),
\end{equation}
where the reduced linear operator is given by
\begin{equation*}
\widetilde{\bA}(\mu)=\Phi_r^T\bA(\mu)\Phi_r =
\begin{pmatrix}
0 & 0 & 0 & 0 \\[0.2cm]
 0 & 0 & f(\mu)  \Phi_{r_u}^{T}\Phi_{r_v} & - \Phi_{r_u}^{T}(D_x\bb)^d\Phi_{r_s} \\[0.2cm]
 0 &- f(\mu) \Phi_{r_v}^{T}\Phi_{r_u} & 0 & - \Phi_{r_v}^{T}(D_y\bb)^d\Phi_{r_s} \\[0.2cm]
0 & 0 & 0 & 0
\end{pmatrix}\in \mathbb{R}^{4r\times 4r},
\end{equation*}
while the reduced quadratic operator is $\widetilde{\bH}=\Phi_r^T\bH(\Phi_r\otimes \Phi_r)\in \mathbb{R}^{4r\times (4r)^2}$. We note that the linear-quadratic structure of the FOM  \eqref{tenfom} is preserved by the reduced system \eqref{tswepfom0} (see \cite{Benner15}) which is also solved in time by Kahan's method \eqref{kahan} as the FOM.

For computational aspects, we separate the quadratic term in \eqref{tswepfom0} and we rewrite the intrusive ROM \eqref{tswepfom0} in the following form
\begin{equation}\label{tswepfom}
\begin{aligned}
\dot{\tilde{\bw}}(t;\mu)= & \; \widetilde{\bA}(\mu)\tilde{\bw}(t;\mu)
+\widetilde{\bH}_1\begin{bmatrix}0\\(\tilde{\bu}(t;\mu)\otimes \tilde{\bu}(t;\mu))\\(\tilde{\bv}(t;\mu)\otimes \tilde{\bu}(t;\mu))\\0\end{bmatrix}
+\widetilde{\bH}_2\begin{bmatrix}(\tilde{\bh}(t;\mu)\otimes \tilde{\bu}(t;\mu))\\(\tilde{\bv}(t;\mu)\otimes \tilde{\bu}(t;\mu))\\(\tilde{\bv}(t;\mu)\otimes \tilde{\bv}(t;\mu))\\(\tilde{\bs}(t;\mu)\otimes \tilde{\bu}(t;\mu)) \end{bmatrix}\\
& +\widetilde{\bH}_3\begin{bmatrix}(\tilde{\bh}(t;\mu)\otimes \tilde{\bv}(t;\mu))\\(\tilde{\bh}(t;\mu)\otimes \tilde{\bs}(t;\mu))\\(\tilde{\bh}(t;\mu)\otimes \tilde{\bs}(t;\mu))\\(\tilde{\bs}(t;\mu)\otimes \tilde{\bv}(t;\mu)) \end{bmatrix}
+\widetilde{\bH}_4\begin{bmatrix}0\\(\tilde{\bs}(t;\mu)\otimes \tilde{\bh}(t;\mu))\\(\tilde{\bs}(t;\mu)\otimes \tilde{\bh}(t;\mu))\\0\end{bmatrix},
\end{aligned}
\end{equation}
where each reduced matrix $\widetilde{\bH}_i\in \mathbb{R}^{4r\times 4r^2}$, $i=1,\ldots ,4$, is in the block diagonal form defined by
\begin{align*}
\widetilde{\bH}_1 &=
\begin{pmatrix}
0 &  &  &  \\
 & \widetilde{\bH}_1^1 &  &   \\
 &  & \widetilde{\bH}_1^2 &   \\
 &  &  &  0
\end{pmatrix} \; , \quad \widetilde{\bH}_2 =
\begin{pmatrix}
\widetilde{\bH}_2^1 &  &  &  \\
 & \widetilde{\bH}_2^2 &  &   \\
 &  & \widetilde{\bH}_2^3 &   \\
 &  &  &  \widetilde{\bH}_2^4
\end{pmatrix}, \\[0.2cm]
\widetilde{\bH}_4 &=
\begin{pmatrix}
0 &  &  &  \\
 & \widetilde{\bH}_4^1 &  &   \\
 &  & \widetilde{\bH}_4^2 &   \\
 &  &  &  0
\end{pmatrix}\; , \quad \widetilde{\bH}_3 = \begin{pmatrix}
\widetilde{\bH}_3^1 &  &  &  \\
 & \widetilde{\bH}_3^2 &  &   \\
 &  & \widetilde{\bH}_3^3 &   \\
 &  &  &  \widetilde{\bH}_3^4
\end{pmatrix},
\end{align*}
with the $r\times r^2$ matrices
\begin{equation} \label{redquad}
\begin{aligned}
\widetilde{\bH}_1^1&=-\Phi_{r_u}^{\top}Q(D_x\Phi_{r_u}\otimes\Phi_{r_u}), \quad \widetilde{\bH}_1^2=-\Phi_{r_v}^{\top}Q(D_x\Phi_{r_v}\otimes\Phi_{r_u}),\\
\widetilde{\bH}_2^1&=-\Phi_{r_h}^{\top}D_xQ(\Phi_{r_h}\otimes\Phi_{r_u}), \quad
\widetilde{\bH}_2^2=-\Phi_{r_u}^{\top}Q(\Phi_{r_v}\otimes D_y\Phi_{r_u}), \\
\widetilde{\bH}_2^3&=-\Phi_{r_v}^{\top}Q(\Phi_{r_v}\otimes D_y\Phi_{r_v}), \quad \widetilde{\bH}_2^4=-\Phi_{r_s}^{\top}Q(D_x\Phi_{r_s}\otimes\Phi_{r_u}),\\
\widetilde{\bH}_3^1&=-\Phi_{r_h}^{\top}D_yQ(\Phi_{r_h}\otimes\Phi_{r_v} ), \quad
\widetilde{\bH}_3^2=-\frac{1}{2} \Phi_{r_u}^{\top}Q(\Phi_{r_h}\otimes D_x\Phi_{r_s}), \\
\widetilde{\bH}_3^3&= -\frac{1}{2} \Phi_{r_v}^{\top}Q(\Phi_{r_h}\otimes D_y\Phi_{r_s}), \quad
\widetilde{\bH}_3^4= -  \Phi_{r_s}^{\top}Q(D_y\Phi_{r_s}\otimes\Phi_{r_v}),\\
\widetilde{\bH}_4^1&=  -\Phi_{r_u}^{\top}Q(\Phi_{r_s}\otimes  D_x\Phi_{r_h}), \quad \widetilde{\bH}_4^2=-\Phi_{r_v}^{\top}Q(\Phi_{r_s}\otimes  D_y\Phi_{r_h}),
\end{aligned}
\end{equation}
where $ Q\in \mathbb{R}^{N\times N^2} $ is the matricized tensor such that $Q (\mathbf{a}\otimes\mathbf{b})=\mathbf{a}\circ\mathbf{b}$ is satisfied for any vectors $\mathbf{a},\mathbf{b}\in\mathbb{R}^{N}$.

All the reduced quadratic operators in \eqref{redquad} are computed in the offline stage. In computing reduced quadratic operators, we follow the tensorial approach in \cite{Karasozen21} by exploiting the structure of Kronecker product for computation of the quadratic operators. For instance, consider the computation of the term $ \Phi_{r_h}^{\top}D_xQ(\Phi_{r_h}\otimes\Phi_{r_u}) $ which can be computed in MATLAB notation as follows
\begin{equation}\label{offline}
\begin{aligned}
 \Phi_{r_h}^{\top}D_xQ(\Phi_{r_h}\otimes\Phi_{r_u})
= \Phi_{r_h}^{\top}D_x
\begin{pmatrix}
\Phi_{r_h}(1,:)\otimes \Phi_{r_u}(1,:)\\
\vdots\\
\Phi_{r_h}(N,:)\otimes \Phi_{r_u}(N,:)
\end{pmatrix}.
\end{aligned}
\end{equation}
Alternatively, the computation of $ G:=Q(\Phi_{r_h}\otimes\Phi_{r_u})\in \mathbb{R}^{N\times r^2}$  can be written as
$$
G(i,:)=\left(\text{vec}\left(\Phi_{r_u}^T \Phi_{r_h}(i,:)\right)\right)^T, \quad i=1,2,\ldots,N,
$$
where vec $ (\cdot) $ denotes the vectorization of a matrix. By using "MULTIPROD" \cite{leva08mmm}, the matrix $G$ can efficiently be computed through the calculation of the tensor $\cG:=\text{MULTIPROD}(\Phi_{r_u},\widehat{\Phi}_{r_h})\in\mathbb{R}^{N\times r \times r}$, where the tensor $\widehat{\Phi}_{r_h}\in \mathbb{R}^{N\times 1\times r}$ is obtained by reshaping the matrix $ \Phi_{r_h} \in \mathbb{R}^{N\times r} $ into $N\times 1\times r$.
Thus, the matrix $G \in \mathbb{R}^{N\times r^2}$ becomes the mode-1 matricization of the tensor $\mathcal{G}$.

\subsection{Operator inference with re-projection }
\label{nonint}

We apply the data-driven, non-intrusive  OpInf method  \cite{Peherstorfer16} with re-projection \cite{Peherstorfer20}  to the RTSWE \eqref{tswe}.
The OpInf  projects trajectories of the  high-dimensional state
spaces of FOMs onto low-dimensional subspaces, and then fit operators to the projected trajectories via least-squares regression. It does not involve any explicit knowledge of the FOM, rather only involves simulation data projected onto the dominant POD subspace. Within OpInf framework, a closure error is introduced into the learned operators, which are different than those of reduced model constructed with the intrusive model reduction. The projected trajectories correspond to non-Markovian dynamics
in the low-dimensional subspaces even though the high dimensional full solutions are Markovian \cite{Chorin06,Stuart04}. In Markovian  dynamical systems 
the future states depend only on the current state and not on previous ones, whereas non-Markovian systems can be considered as dynamical systems with a memory so that future states depend on the current and previous states. In \cite{Peherstorfer20}, a data sampling scheme is introduced that cancels these non-Markovian dynamics after each time step such that the trajectories correspond to Markovian dynamics in the reduced space. In this way,  the OpInf exactly recovers the reduced model from these re-projected trajectories. The data sampling scheme in \cite{Peherstorfer20} iterates between time stepping the high-dimensional FOM and projecting onto low-dimensional reduced subspaces to retain
the low-dimensional Markovian dynamics. It was shown in \cite{Peherstorfer20} that under certain conditions, applying OpInf to re-projected trajectories gives the same operators that are obtained with the intrusive model reduction in the limit of $r \rightarrow  N$ for many dynamical systems with polynomial nonlinear terms. This is validated in the numerical results in Section \eqref{sec:num} so that the learnt low-dimensional models have almost the same level of accuracy as the intrusive reduced models.

Using the definitions form the previous section, the reduced trajectories for each state are denoted as
\begin{align*}
&\widehat{W}_h(\mu)=\Phi_{r_h}^{\top} W_h(\mu) \in\mathbb{R}^{r\times K}, \quad \widehat{W}_u(\mu)=\Phi_{r_u}^{\top} W_u(\mu) \in\mathbb{R}^{r\times K},\\
&\widehat{W}_v(\mu)=\Phi_{r_v}^{\top} W_v(\mu) \in\mathbb{R}^{r\times K}, \quad \widehat{W}_s(\mu)=\Phi_{r_s}^{\top} W_s(\mu) \in\mathbb{R}^{r\times K}.
\end{align*}
Furthermore, the time derivatives and their projections are given by
\begin{align*}
&\dot{W_h}(\mu) =  \left[\dot{\bh}^1(\mu),\ldots,\dot{\bh}^K(\mu)\right], \quad \dot{W_u}(\mu)  =  \left[\dot{\bu}^1(\mu),\ldots,\dot{\bu}^K(\mu)\right],\\
&\dot{W_v}(\mu) =  \left[\dot{\bv}^1(\mu),\ldots,\dot{\bv}^K(\mu)\right], \quad \dot{W_s}(\mu) =  \left[\dot{\bs}^1(\mu),\ldots,\dot{\bs}^K(\mu)\right],\\
&\dot{\widehat{W}}_h(\mu)=\Phi_{r_h}^{\top}\dot{ W_u}(\mu) , \quad \dot{\widehat{W}}_u(\mu)=\Phi_{r_u}^{\top}\dot{ W_u}(\mu) ,\\
&\dot{\widehat{W}}_v(\mu)=\Phi_{r_v}^{\top} \dot{W_v}(\mu) , \quad \dot{\widehat{W}}_s(\mu)=\Phi_{r_s}^{\top} \dot{W_s}(\mu).
\end{align*}

We compute the time derivatives by evaluating the right-hand side of the RTSWE \eqref{tenfom}. Alternatively, the time derivatives can be approximated using finite differences \cite{Martins13}.
The reduced operators are recovered by solving the following least-squares problems separately for each state vector
\begin{equation}\label{eq:lsqr}
\min_{\cX_{j}}\sum_{k = 1}^M\left\|\cA_{j}(\mu_k)\cX_{j}^{\top} - \dot{\widehat{W}}_{j}(\mu_k)^{\top}\right\|_{F}^{2},\quad j=h,u,v,s,
\end{equation}
where the data matrices are defined by
\begin{align*}
&\cA_{h}(\mu_k)=\left[(\widehat{W}_h(\mu_k)\hat{\otimes}\widehat{W}_u(\mu_k))^{\top},(\widehat{W}_h(\mu_k)\hat{\otimes}\widehat{W}_v(\mu_k))^{\top}\right]\in \mathbb{R}^{K\times 2r^2}, \\
&\cA_{u}(\mu_k)=\left[(\widehat{W}_u(\mu_k)\hat{\otimes}\widehat{W}_u(\mu_k))^{\top},(\widehat{W}_v(\mu_k)\hat{\otimes}\widehat{W}_u(\mu_k))^{\top},(\widehat{W}_h(\mu_k)\hat{\otimes}\widehat{W}_s(\mu_k))^{\top},f(\mu_k)\widehat{W}_v(\mu_k)^{\top}\right]\in\mathbb{R}^{K\times (3r^2+r)},\\
&\cA_{s}(\mu_k)=\left[(\widehat{W}_u(\mu_k)\hat{\otimes}\widehat{W}_s(\mu_k))^{\top},(\widehat{W}_v(\mu_k)\hat{\otimes}\widehat{W}_s(\mu_k))^{\top}\right]\in \mathbb{R}^{K\times 2r^2}, \\
&\cA_{v}(\mu_k)=\left[(\widehat{W}_u(\mu_k)\hat{\otimes}\widehat{W}_v(\mu_k))^{\top},(\widehat{W}_v(\mu_k)\hat{\otimes}\widehat{W}_v(\mu_k))^{\top},(\widehat{W}_h(\mu_k)\hat{\otimes}\widehat{W}_s(\mu_k))^{\top},f(\mu_k)\widehat{W}_u(\mu_k)^{\top}\right]\in\mathbb{R}^{K\times (3r^2+r)},
\end{align*}
where $\hat{\otimes}$ stands for the column-wise Kronecker product (Khatri-Rao product), and the learned operators $\cX_{j}$ have the following form
\begin{align*}
\cX_{h}=\left[\widehat{\cH}_{1,1},\widehat{\cH}_{1,2}\right]\in \mathbb{R}^{r\times 2r^2}, \quad \cX_{u}=\left[\widehat{\cH}_{2,1},\widehat{\cH}_{2,2},\widehat{\cH}_{2,3}, \widehat{\cL}_{1}(\mu)\right]\in \mathbb{R}^{r\times (3r^2+r)},\\
\cX_{s}=\left[\widehat{\cH}_{4,1},\widehat{\cH}_{4,2}\right]\in \mathbb{R}^{r\times 2r^2}, \quad \cX_{v}=\left[\widehat{\cH}_{3,1},\widehat{\cH}_{3,2},\widehat{\cH}_{3,3}, \widehat{\cL}_{2}(\mu)\right]\in \mathbb{R}^{r\times (3r^2+r)}.
\end{align*}
The least-squares problem \eqref{eq:lsqr} can be written in compact form as
\begin{equation}\label{eq:clsqr}
\min_{\cX_{j}}\left\|\cA_{\mu j}\cX_{j}^{\top} - \dot{\widehat{W}}_{\mu j}\right\|_{F}^{2},\quad j=h,u,v,s,
\end{equation}
where $ \dot{\widehat{W}}_{\mu j}=\left[\dot{\widehat{W}}_j^T(\mu_1),\ldots,\dot{\widehat{W}}_j^T(\mu_M)\right]^T  $ and $\cA_{\mu j}=\left[\cA_{j}^T(\mu_1),\ldots,\cA_{j}^T(\mu_M)\right]^T$ for $ j=h,u,v,s $.

We remark that each of the least squares problems \eqref{eq:lsqr} are independent and are solved separately to increase the computational efficiency in the online stage \cite{Peherstorfer16}.

The ROM in the OpInf framework is given as
\begin{align*}	
	\dot{\widehat{\bh}}(t;\mu) &  =  \widehat{\cH}_{1,1}(\widehat{\bh}(t;\mu)\otimes\widehat{\bu}(t;\mu))  +\widehat{\cH}_{1,2}(\widehat{\bh}(t;\mu)\otimes \widehat{\bv}(t;\mu))   \\
	\dot{\widehat{\bu}}(t;\mu) & =  \widehat{\cH}_{2,1}(\widehat{\bu}(t;\mu)\otimes\widehat{\bu}(t;\mu))  +\widehat{\cH}_{2,2}(\widehat{\bv}(t;\mu)\otimes\widehat{\bu}(t;\mu)) +\widehat{\cH}_{2,3}(\widehat{\bh}(t;\mu)\otimes\widehat{\bs}(t;\mu)) + f(\mu)\widehat{\cL}_{1}\widehat{\bv}(t;\mu) \\
	\dot{\widehat{\bv}}(t;\mu) & =  \widehat{\cH}_{3,1}(\widehat{\bu}(t;\mu)\otimes\widehat{\bv}(t;\mu))  +\widehat{\cH}_{3,2}(\widehat{\bv}(t;\mu)\otimes\widehat{\bv}(t;\mu)) +\widehat{\cH}_{3,3}(\widehat{\bh}(t;\mu)\otimes\widehat{\bs}(t;\mu)) + f(\mu) \widehat{\cL}_{2} \widehat{\bu}(t;\mu) \\
	\dot{\widehat{\bs}}(t;\mu) & =   \widehat{\cH}_{4,1}(\widehat{\bu}(t;\mu)\otimes\widehat{\bs}(t;\mu)) +\widehat{\cH}_{4,2}(\widehat{\bv}(t;\mu)\otimes\widehat{\bs}(t;\mu)).
\end{align*}

Data sampling via re-projection \cite{Peherstorfer20}
is summarized in Algorithm \eqref{alg:ReProj} by maintaining the coupling structure of the ROM as in \cite{Kramer20}, while the process of the operator inference with re-projection for parametric RTSWE is given in Algorithm \eqref{algo1}.

\begin{algorithm}[htb!]
	\caption{Data sampling via re-projection \label{alg:ReProj}}	
	\begin{algorithmic}[1]
	\Procedure{Re-Projection}{}\\
	\hspace*{\algorithmicindent} \textbf{Input:} States $ \bw^j $ $ j=1,\ldots,K $ and POD basis matrix $ \Phi_r $
    \For{$j=1,\ldots,K$}{\\
    \qquad \quad Set $\bw^j_{\text{proj}}=(\bh^j_{\text{proj}},\bu^j_{\text{proj}}, \bv^j_{\text{proj}}, \bs^j_{\text{proj}})^T=\Phi_r \Phi_r^T  \bw^j $\\
    \qquad \quad  Query FOM at each time step for $ \dot{\overline{\bw}}^j=\Phi_r^T\bF(\bw^j_{\text{proj}}) $ in \eqref{tenfom}\\
    \qquad \quad Re-project the states $ \overline{\bw}^{j}= \Phi_r^T \bw^{j}_{\text{proj}}$
    }
    \EndFor\\
    \Return $ \dot{\overline{W}}=\left[\dot{\overline{\bw}}^1,\ldots,\dot{\overline{\bw}}^K\right]  $ and $ \overline{W}=\left[\overline{\bw}^1,\ldots,\overline{\bw}^K\right] $
    \EndProcedure
	\end{algorithmic}
\end{algorithm}

\begin{algorithm}[htb!]
	\caption{Operator inference with re-projection for parametric RTSWE \label{algo1}}
	\begin{algorithmic}[1]
		\Procedure{OpInf}{}\\
		 \hspace*{\algorithmicindent} \textbf{Input:} States $ \bw_j(\mu_i) $ and time derivatives $\dot{\bw}_j(\mu_i) $ for  $i = 1,\ldots,M$ and $ j=1,\ldots,K $
		\State Construct the trajectories for each state
		\begin{align*}
		&W_h(\mu) =  \left[\bh_1(\mu),\ldots,\bh_K(\mu)\right], \quad W_u(\mu)  =  \left[\bu_1(\mu),\ldots,\bu_K(\mu)\right],\\
		&W_v(\mu) =  \left[\bv_1(\mu),\ldots,\bv_K(\mu)\right], \quad W_s(\mu) =  \left[\bs_1(\mu),\ldots,\bs_K(\mu)\right].
		\end{align*}
		\State Construct the global snapshot matrices  $W_{j\mu} =  \left[W_j(\mu_1),\ldots,W_j(\mu_M)\right] $ for $ j=h,u,v,s $
		\State Compute the global POD basis $ \Phi_{r_j} $ of $W_{j\mu}$ for $ j=h,u,v,s $
		\State Sample the data via re-projection   and set $ \widehat{W}_j=\overline{W}_j $ and $ \dot{\widehat{W}}_j=\dot{\overline{W}}_j $ for $ j=h,u,v,s $
		\State Determine the tolerance of \texttt{lsqminnorm} $ tol $, using L-curve formula
     \State Solve the least-squares problem \eqref{eq:lsqr} to obtain operators of the reduced-order system using $ \cX_{j}^T =$\texttt{lsqminnorm}$ (\cA_{\mu j},\dot{\widehat{W}}_{\mu j},tol) $
		\EndProcedure
	\end{algorithmic}
\end{algorithm}

For quadratic systems like the RTSWE, the number of the time steps should satisfy $K \ge r + r^2$, so that the least square system  \eqref{eq:lsqr} is overdetermined and a unique solution exits according to the Corollary 3.2 in \cite{Peherstorfer20} and the data matrices  $\cA_{j}, \;j=h,u,v,s$ have full rank.
When this condition is not fulfilled, the least-squares problem \eqref{eq:clsqr}  is underdetermined and therefore uniqueness of the solution is lost.
The least-squares problem arising in OpInf can be ill-conditioned; thus, regularization should be applied. The common regularization techniques are the truncated SVD and the Tikhanov regularization \cite{tikhonov77}. An ill-conditioned matrix as either a matrix with a well-determined numerical rank or an ill-determined numerical rank, depending on the behavior of the singular value spectrum  \cite{Hansen87}. Usually, the L-curve criterion \cite{Hansen00} is used to determine the tolerance at which the singular values are truncated.  Analyzing the L-curve, a tolerance can be determined that has a good compromise between the matching of data-fidelity term and making the problem well-conditioned.
As it will be illustrated in Section \ref{sec:num}, the singular values of the data matrices decay monotonically without exhibiting a gap.
According to terminology in \cite{Hansen87}, the numerical rank of the data matrices is ill-determined. The L-curve also does not provide any information; there exist not an optimal value for the tolerance, consequently truncated SVD can not be used.
 One way to overcome this bottleneck is to solve the least-squares problem via the minimum norm of the solution. The minimum-norm solution enforces uniqueness of the solution by picking smallest $ \|\cX_{j}\|_{F}  $ which minimizes $ \|\cA_{\mu j}\cX_{j}^{\top} - \dot{\widehat{W}}_{\mu j}\|_{F}$. To obtain the minimum-norm solution, one can use the Moore-Penrose pseudoinverse or complete orthogonal decomposition (COD) formulas. We use here MATLAB's routine {\texttt lsqminnorm} as regularizer of the least-squares problem \eqref{eq:lsqr} to achieve the uniqueness. The tolerance  will be determined with the L-curve as shown in Section \ref{sec:num}.

Under the conditions that the time stepping scheme for the FOM is convergent as $\Delta t \longrightarrow 0$, the derivatives approximated from the projected states $\dot{\hat{\bw}}$ converge to $\frac{d}{dt}\hat{\bw}(t_k)$ as $\Delta t \longrightarrow 0$, and the data matrix is of full rank, then the  learned operators  $\hat{A},\hat{H}$ converge to the intrusively projected operators $\tilde{A},\tilde{H}$ \cite{Kramer20,Peherstorfer16}. Kahan's method is second order convergent, the data matrix is of full rank, therefore these conditions above are fulfilled.

In the OpInf framework, once the operators are inferred, any time stepping scheme can be used to solve the ROM with the inferred operators, because  the operators of the system of ODEs are inferred and not the operators of the time-discretized system. Here, we use again Kahan's method as time-discretization scheme for \eqref{kahan}, because it is cheap, i.e., linearly implicit, and preserves the conserved quantities in contrast to the implicit or explicit Euler methods.

\subsection{Computational aspects}

The costs of the data generation and the construction of the POD basis are the same for the POD-G and the OpInf. We, therefore, consider the costs of the  OpInf Algorithm \eqref{algo1} and compare it to the cost of the intrusive MOR in Section \eqref{POD-G}.

The computational cost of the OpInf  Algorithm \eqref{algo1} is dominated by querying the high-dimensional system. The costs of the projection of the trajectories onto the POD space are of order $\mathcal{O}(rNK)$ and thus are bounded linearly in the dimension $N$ of the FOM and the number of time steps $K$ for each parameter value $ \mu_k $.
The data matrices $\cA_{j}, \;j=h,u,v,s$ and the right-hand side matrices  $\dot{\widehat{W}}_{j}, \; j=h,u,v,s$  are assembled with the costs at most $\mathcal{O}(MK( 3r^2+r)$) and $\mathcal{O}(MKr)$, respectively.
The leading costs of solving one of the $r$ independent least-squares problems for each state vector $h,u,v,s$  is $\mathcal{O}(KM^6r^6)$, see, e.g., \cite{Peherstorfer16}.
Overall, the cost of Algorithm   \eqref{algo1} is bounded linearly in the dimension $N$ of the FOM and linearly in the number of time steps $K$. Only the projection of the trajectories onto the POD space incurs cost that scale with the dimension $N$  of the FOM. All other computations in Algorithm \eqref{algo1} have cost independent of the dimension $N$ of the FOM.  Additional cost occurs in generating the re-projected trajectories and the time derivatives, which depends on the number of time steps $K$ used for generating the re-projected trajectories. Therefore, the computational cost of OpInf Algorithm \eqref{algo1} with re-projection is twice of OpInf without re-projection.

A data matrix with a large condition number can introduce significant numerical errors into the solutions of the $r$ least-squares problems, and thus into the inferred operators. Increasing the dimension $r$ of the reduced space decreases the ROM error,  but increases the degrees of freedom in OpInf, so that the data matrix becomes ill-conditioned.
A large condition number typically arises if the states
at different time steps are almost linearly dependent in the case of parametric systems. There exist several heuristic strategies that might help to reduce the condition number of the data matrix.
The condition number of the least-squares problem can be decreased by sampling in a subset of trajectories as well as the computational efficiency can be increased, by taking every $k$-th snapshot in the data matrix \cite{Peherstorfer16}. The condition number decreases adding more training data,  starting to solve the FOM with different initial values \cite{Peherstorfer16} because the first $r$ POD basis becomes richer \cite{swischuk2019projection}. Both strategies are employed in Section \ref{sec:num}.

\section{Numerical results}
\label{sec:num}

In this section, we compare the accuracy, computational efficiency and the prediction capabilities of the POD-G and OpInf in the parametric and in the non-parametric form for the double vortex RTSWE \cite{Eldred19} in a doubly periodic domain $\Omega =[0, L]^2$ without the bottom topography $(b =0)$.
The initial conditions are given by
\begin{align*}
&h_0(x,y)=H_0-\Delta h\left[e^{-0.5((x'_1)^2+(y'_1)^2)}+e^{-0.5((x'_2)^2+(y'_2)^2)}-\frac{4\pi \sigma_x \sigma_y}{L^2}\right],\\
&u_0(x,y)=\frac{-g\Delta h}{f \sigma_y}\left[y_1''e^{-0.5((x'_1)^2+(y'_1)^2)}+y''_2e^{-0.5((x'_2)^2+(y'_2)^2)}\right],
\\
&v_0(x,y)=\frac{g\Delta h}{f \sigma_x}\left[x_1''e^{-0.5((x'_1)^2+(y'_1)^2)}+x''_2e^{-0.5((x'_2)^2+(y'_2)^2)}\right],
\\
&s_0(x,y)=g\left(1+0.05\sin\left[\frac{2\pi}{L}(x-xc)\right]\right),
\end{align*}
where $ xc=0.5L $ and
\begin{align*}
	&x'_1=\frac{L}{\pi \sigma_x}\sin\left[\frac{\pi}{L}(x-xc_1)\right],\qquad x'_2=\frac{L}{\pi \sigma_x}\sin\left[\frac{\pi}{L}(x-xc_2)\right],\\
	&y'_1=\frac{L}{\pi \sigma_y}\sin\left[\frac{\pi}{L}(y-yc_1)\right],\qquad y'_2=\frac{L}{\pi \sigma_y}\sin\left[\frac{\pi}{L}(y-yc_2)\right]\\
	&x''_1=\frac{L}{2\pi \sigma_x}\sin\left[\frac{2\pi}{L}(x-xc_1)\right],\qquad x''_2=\frac{L}{2\pi \sigma_x}\sin\left[\frac{2\pi}{L}(x-xc_2)\right],\\
	&y''_1=\frac{L}{2\pi \sigma_y}\sin\left[\frac{2\pi}{L}(y-yc_1)\right],\qquad y''_2=\frac{L}{2\pi \sigma_y}\sin\left[\frac{2\pi}{L}(y-yc_2)\right].
\end{align*}
The center of the two vortices are given by
$$
xc_1=(0.5-ox)L,\; xc_2=(0.5+ox)L,\; yc_1=(0.5-oy)L,\; yc_2=(0.5+oy)L.
$$
We set parameter values $ L=5000  $km,  $ H_0=750 $m,
$ h =75 $m, $ g=9.80616 $ms$ ^{-2} $, $ \sigma_x=\sigma_y= 3L/40$ and $ ox =oy =0.1 $.
All simulations are performed on a grid with the mesh sizes  $n_x  = n_y  = 120$km on a machine with Intel$^{\circledR}$
Core$^{{\mathrm TM}}$ i7 2.5 GHz 64 bit CPU, 8 GB RAM, Windows 10, using 64 bit MatLab R2019a.

\subsection{The non-parametric case}

For the non-parametric case, we take Coriolis parameter as  $f=0.00006147\text{s}^{-1}$ at the latitude $7$, close to equator and the number of time steps are set to $K=250$ with the time step-size $\Delta t= 486$s, that leads to the final time $T=33$h $45$min \cite{Eldred19}.
The size of  each snapshot matrix of  the FOM is $14440 \times 250$.

The accuracy of the ROMs is measured using the relative error of the four states
\begin{align}\label{eq:rel_err}
\frac{||\Phi_r Z-W||_{F}}{||W||_{F}},
\end{align}
where $ W\in\mathbb{R}^{4N\times K} $ is the snapshot matrix of the FOM and  $Z\in\mathbb{R}^{4r\times K} $ is the snapshot matrix of either the non-intrusive or the intrusive ROMs.
The average relative errors between FOM and ROM solutions are given for each state variable  ${\bm w}= {\bm u}, {\bm v},{\bm h}, {\bm s} $ in the  time-averaged $L^2(\Omega)$-norms as
\begin{align}\label{relerrorsol}
\|\bm{w}-{\bm   w}_r\|_{\mathrm rel}=\frac{1}{K}\sum_{k=1}^{N_t}\frac{\|{\bm  w}^k-{\bm  w}_r^k\|_{L^2}}{\|{\bm  w}^k\|_{L^2}}, \quad  \|{\bm  w}^k\|_{L^2}^2=\sum_{i=1}^n \sum_{j=1}^n(w_{i,j}^k)^2\Delta x\Delta y,
\end{align}
where ${\bm  w}_r$ denotes the reduced approximation to the full  solution  ${\bm w}$, either of the G-POD $\widehat{\bm  w}$ or OpInf  $\widetilde{\bm  w}$.

Conservation of the discrete conserved quantities \eqref{dener}: the energy, mass, buoyancy, and total vorticity obtained by the FOM and ROM solutions are measured using the time-averaged relative error defined by
\begin{align}\label{abserrorcon}
\frac{1}{K}\sum_{k=1}^{K}\frac{|E(\bm{w}^k)-E(\bm{w}^0)|}{|E(\bm{w}^0)|}, \quad E\equiv H,M,B,Q.
\end{align}

\begin{figure}[htb!]
	\centering
	\includegraphics[width=0.45\columnwidth]{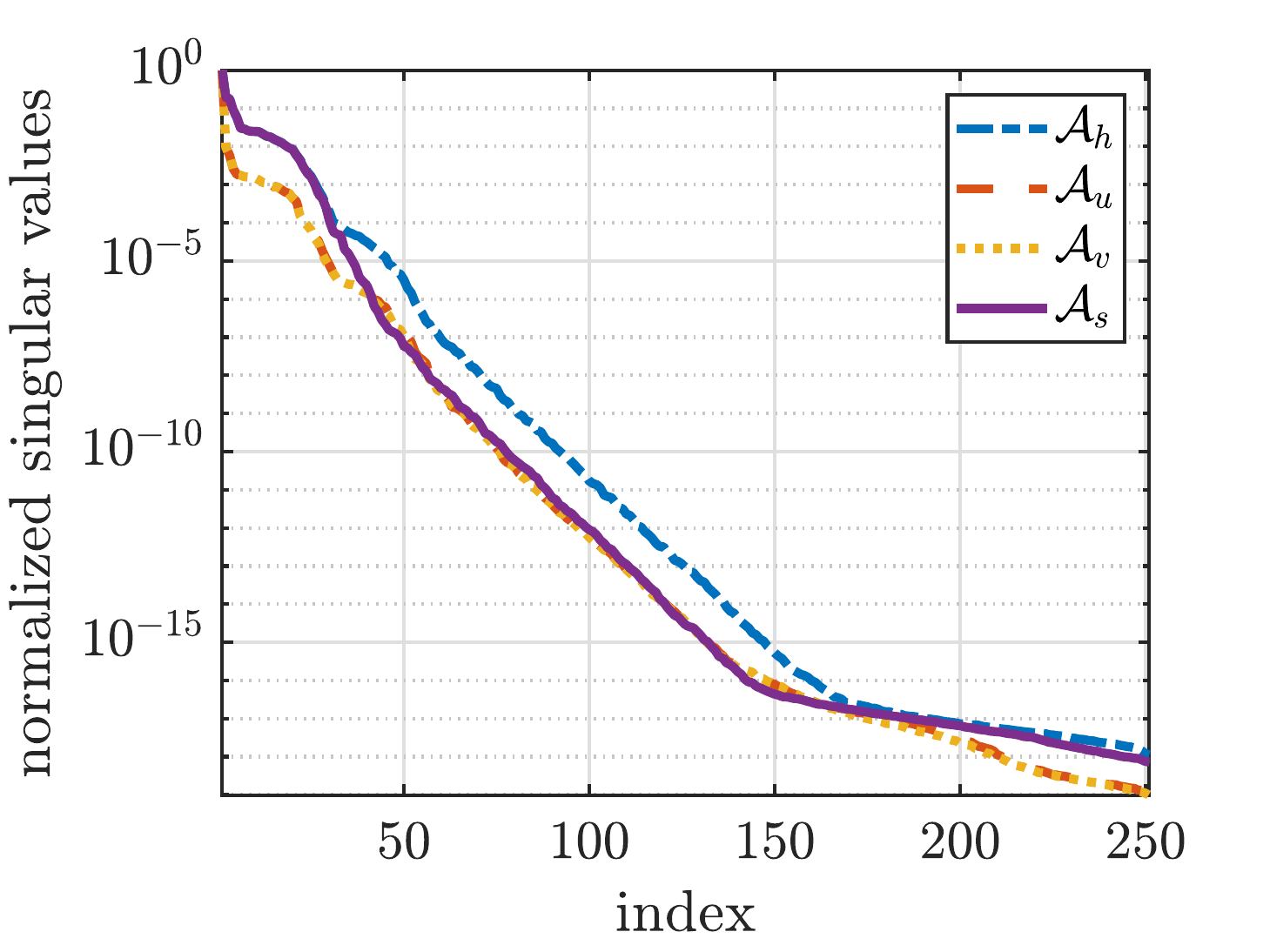}
	\caption{Normalized singular values of the data matrices for the reduced dimension  $r=20$. \label{fig:datasing}}	
\end{figure}

Firstly, we illustrate  how the tolerance $tol$ in Algorithm~\ref{algo1} is determined for the regularization of the least-squares problem \eqref{eq:clsqr}.
Singular values of the data matrices in Figure \ref{fig:datasing} decay without exhibiting a gap. Therefore, the L-curve does not provide any information to determine the tolerances for the truncated SVD.

\begin{figure}[H]
	\centering
	\includegraphics[width=.45\columnwidth]{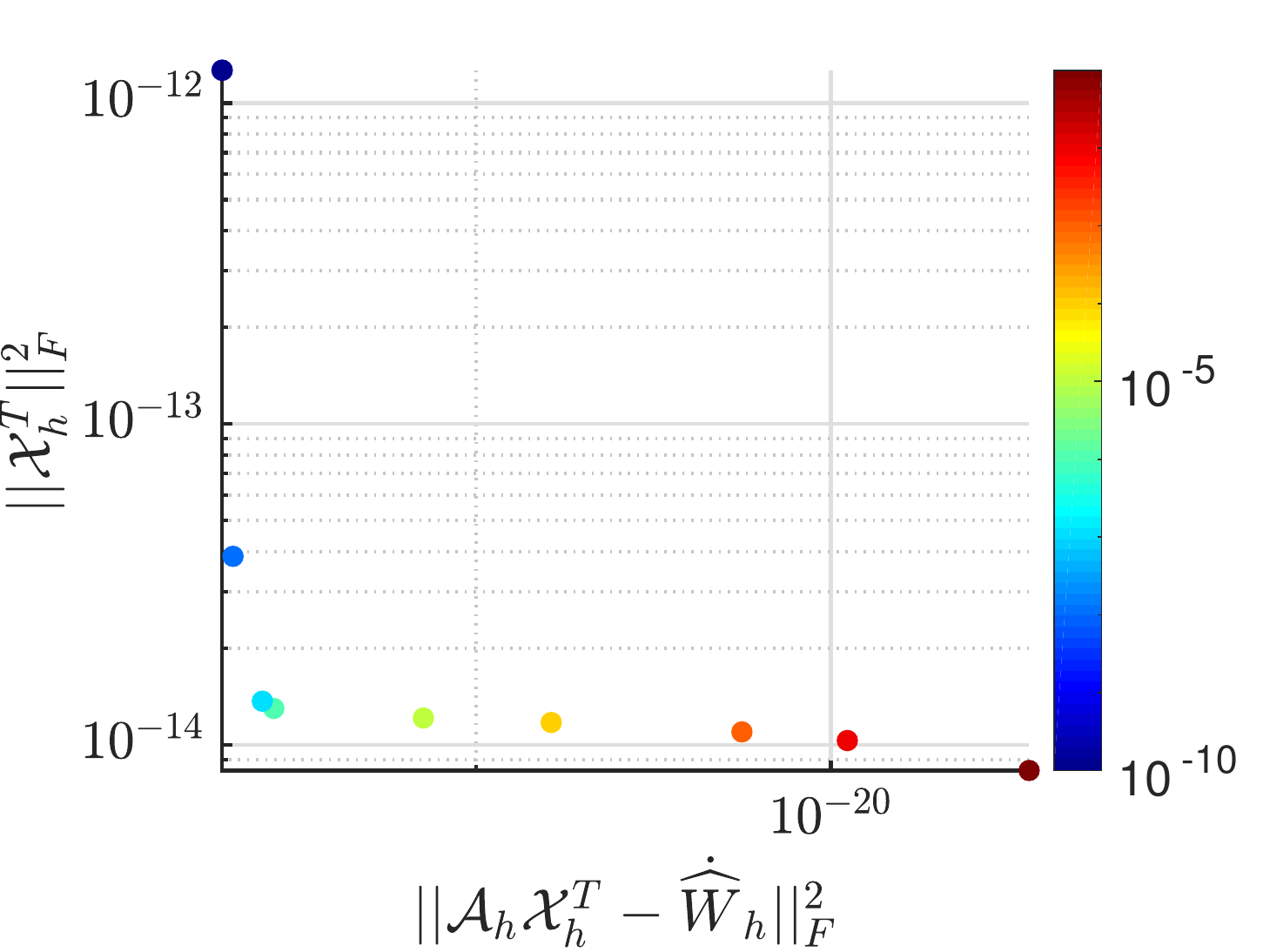}
   \includegraphics[width=.45\columnwidth]{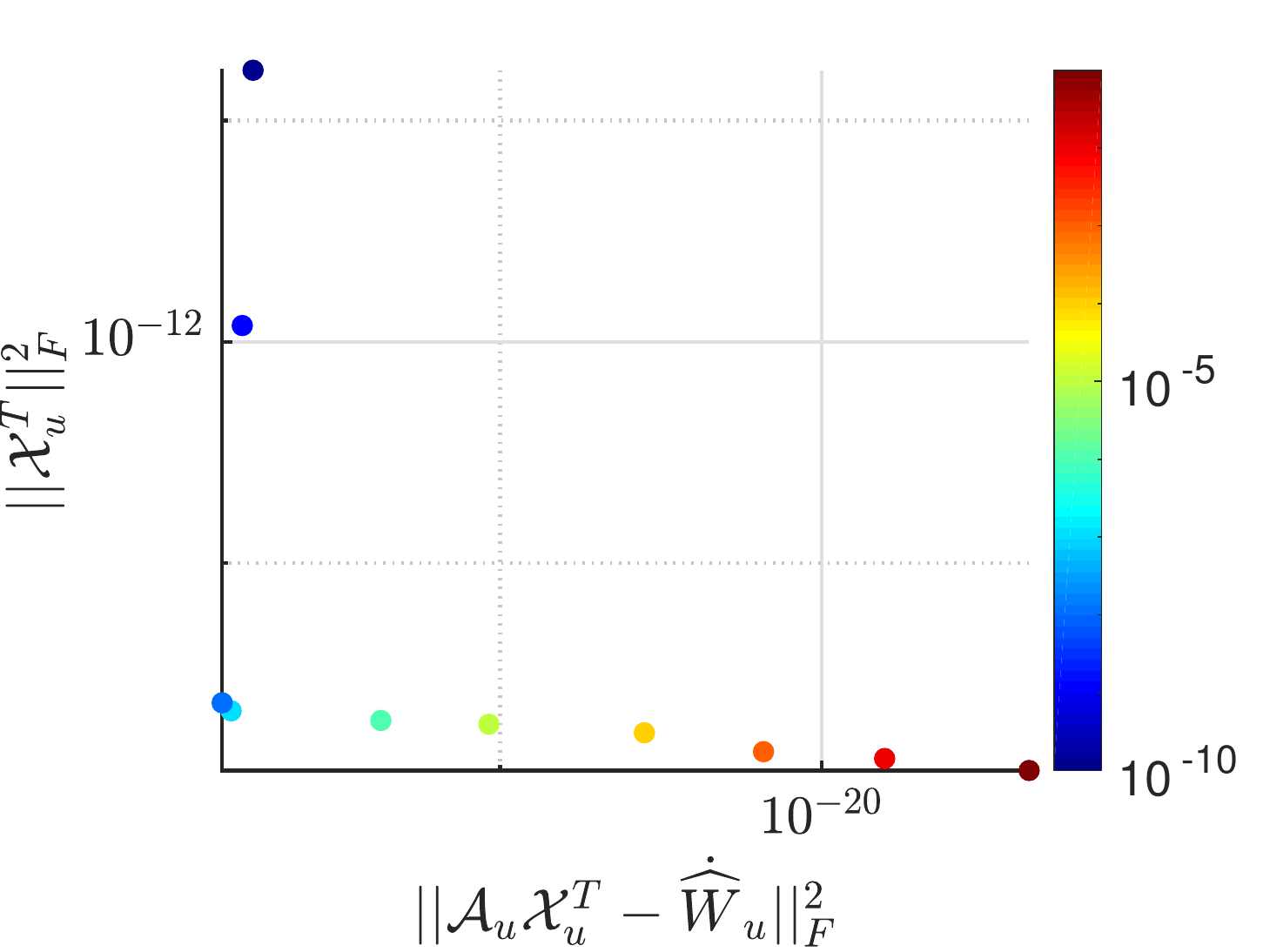}
	\includegraphics[width=.45\columnwidth]{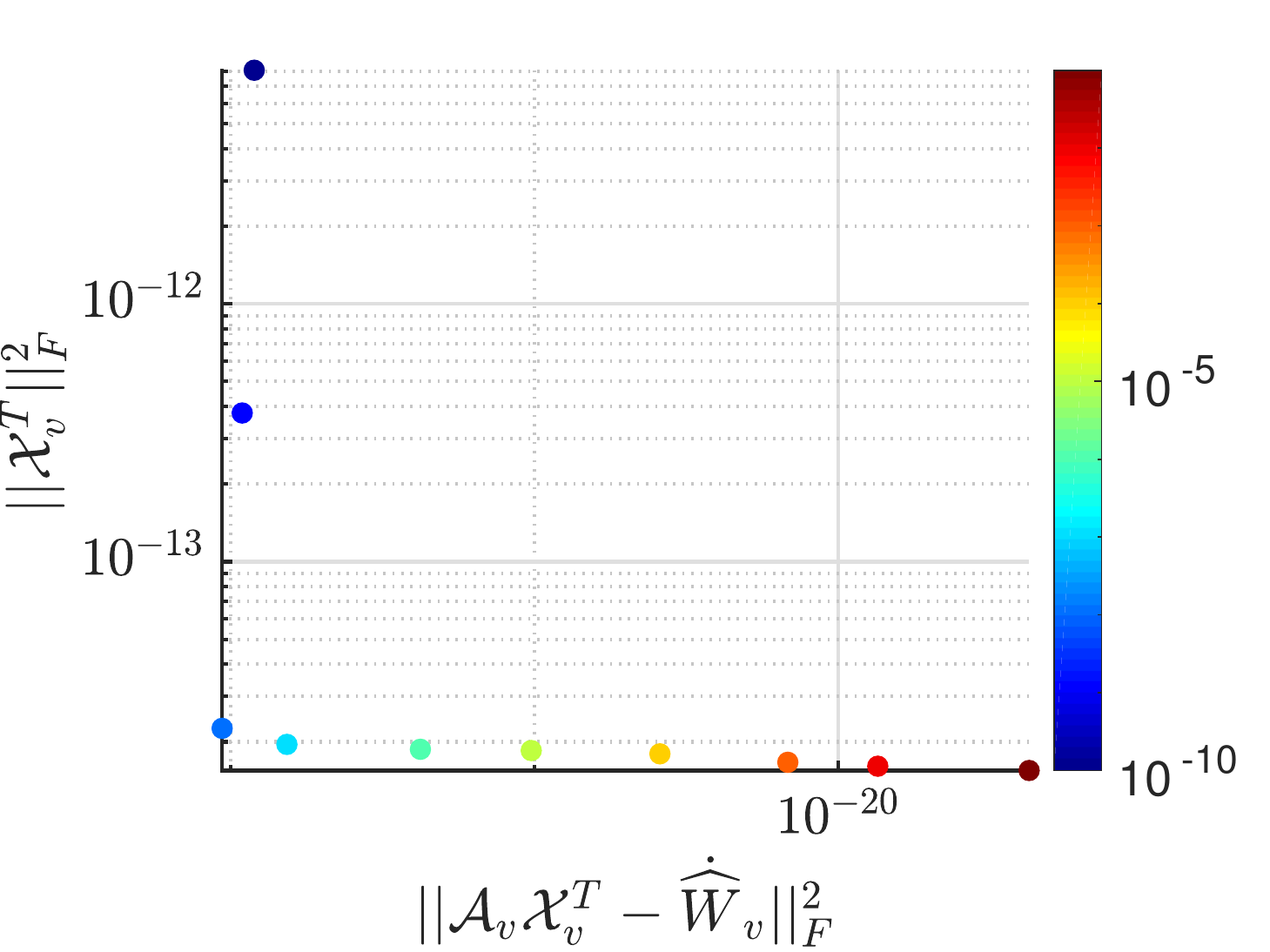}	
	\includegraphics[width=.45\columnwidth]{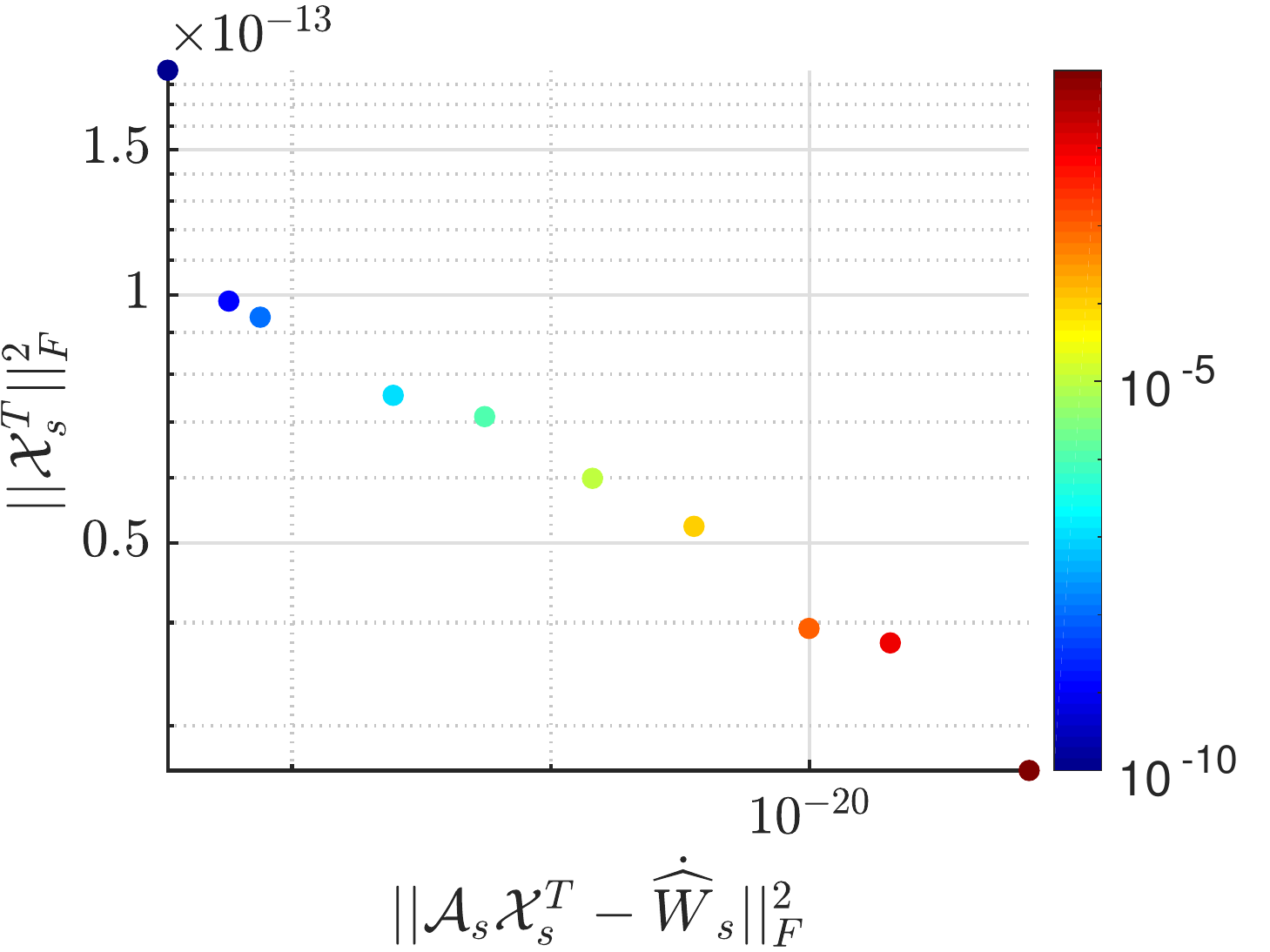}
	\caption{L-curves for reduced dimension $ r=20. $ \label{fig:lcurve}}	
\end{figure}

The tolerances of the least-squares minimum norm solver routine {\texttt lsqminnorm } are determined by the L-curves of the data matrices in Figure~\ref{fig:lcurve} as $ 1e-6 $ for $ j=h,v $ and $ 1e-7 $ for $ j=u,s$.

\begin{figure}[H]
	\centering
	\includegraphics[width=0.45\columnwidth]{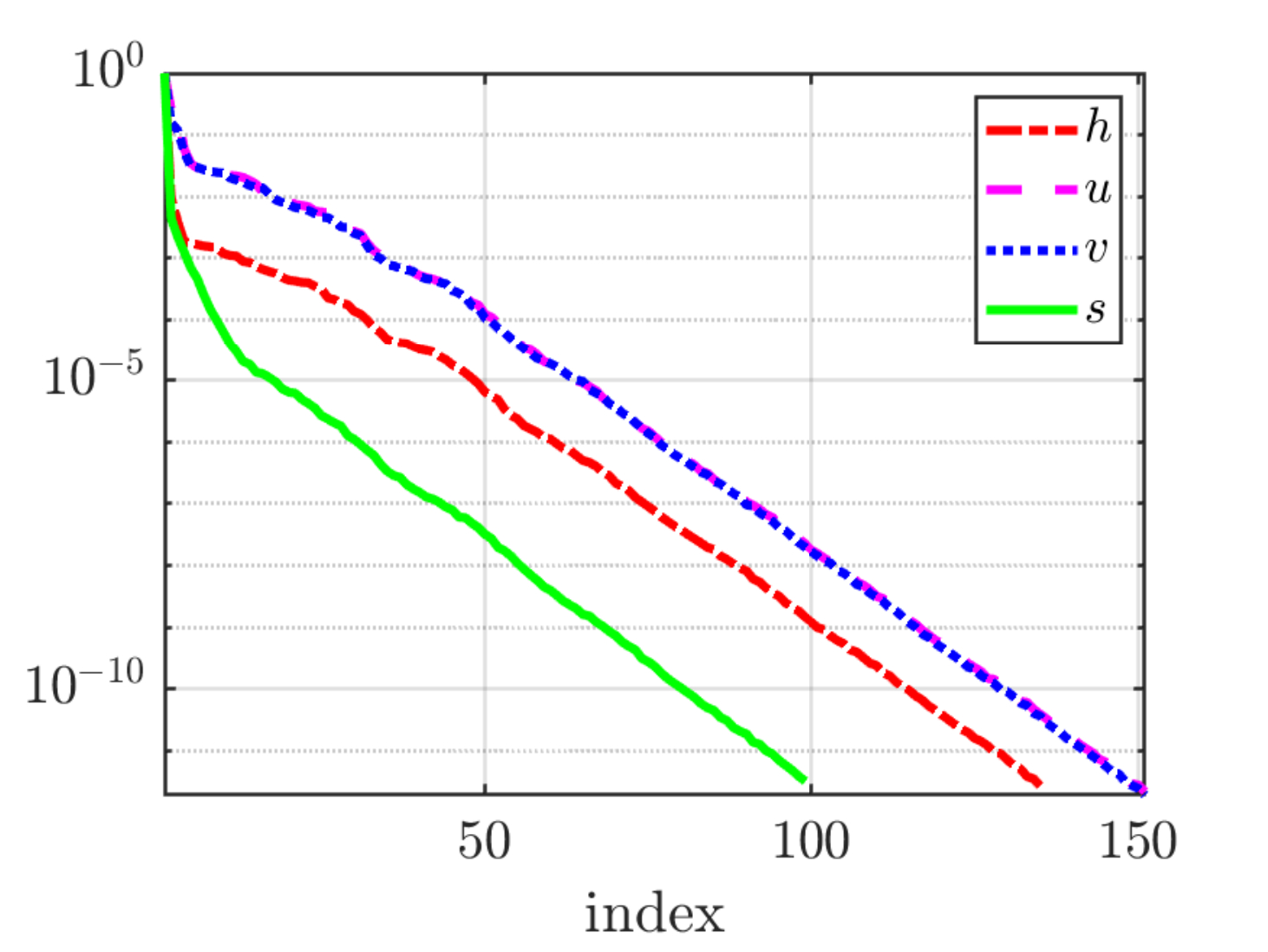}
	\caption{Normalized singular values of the snapshot matrices. \label{fig:tswesing}}	
\end{figure}

In Figure~\ref{fig:tswesing}, the singular values decay relatively  slowly, as for the complex fluid dynamic problems with wave and transport phenomena like the SWEs \cite{Ohlberger16}.

\begin{figure}[H]
	\centering
	\includegraphics[width=0.45\columnwidth]{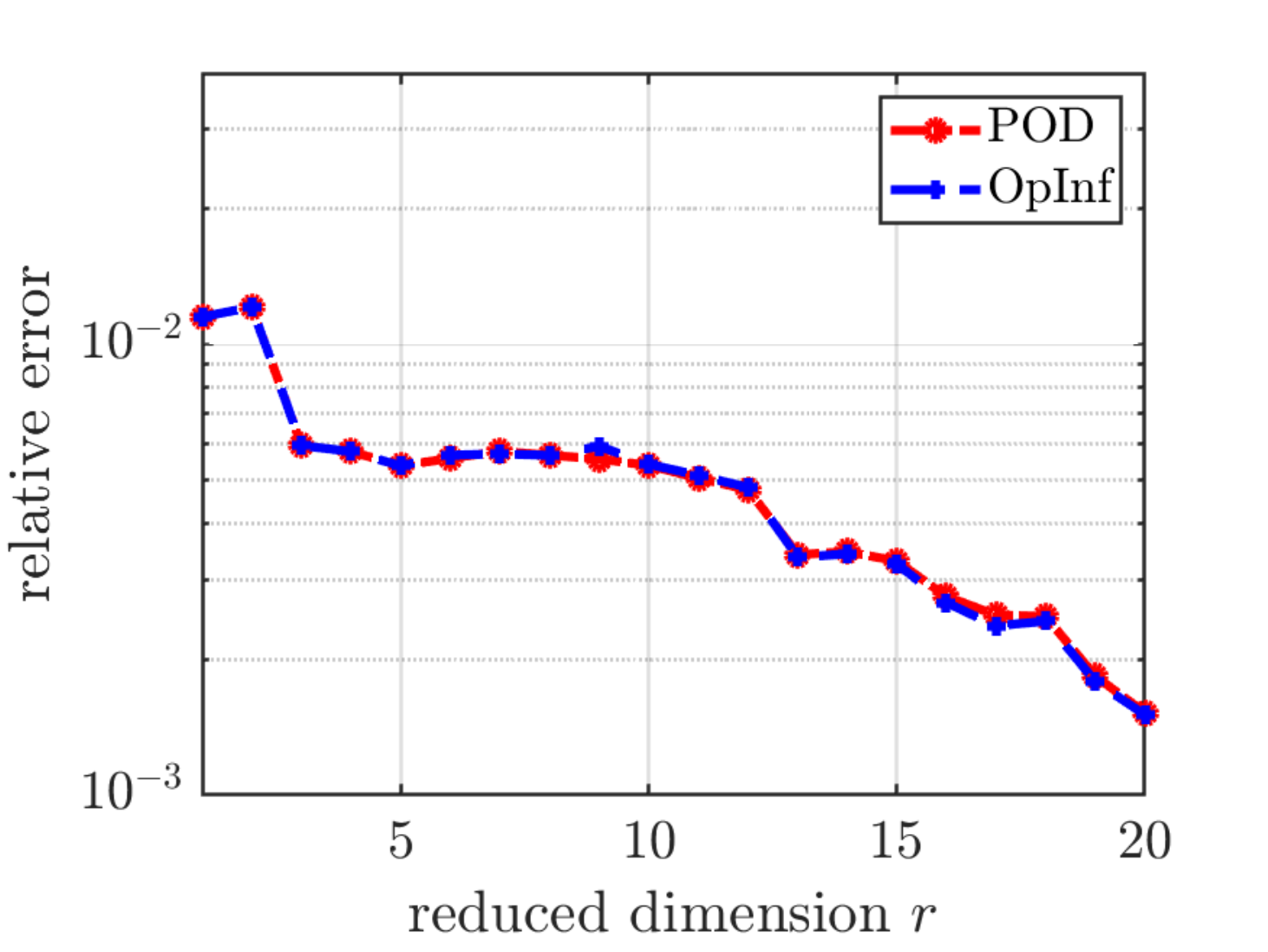}
	\caption{Relative errors \eqref{eq:rel_err} of the POD-G and OpInf. \label{fig:relerr}}	
\end{figure}

In  Figure \ref{fig:relerr}, the relative FOM-ROM errors \eqref{eq:rel_err} of the four states are plotted. The errors of the  POD-G and OpInf are very close.
With increasing reduced dimension $r$, the POD-G and OpInf errors decrease. Average relative errors \eqref{relerrorsol} of the POD-G  and OpInf are $1.499e-03 $ and $1.485e-03 $, relatively, for the reduced dimension $r=20$, which shows that both ROMs approximate the FOM with the same level of accuracy.

\begin{figure}[H]
	\centering
	\includegraphics[width=0.45\columnwidth]{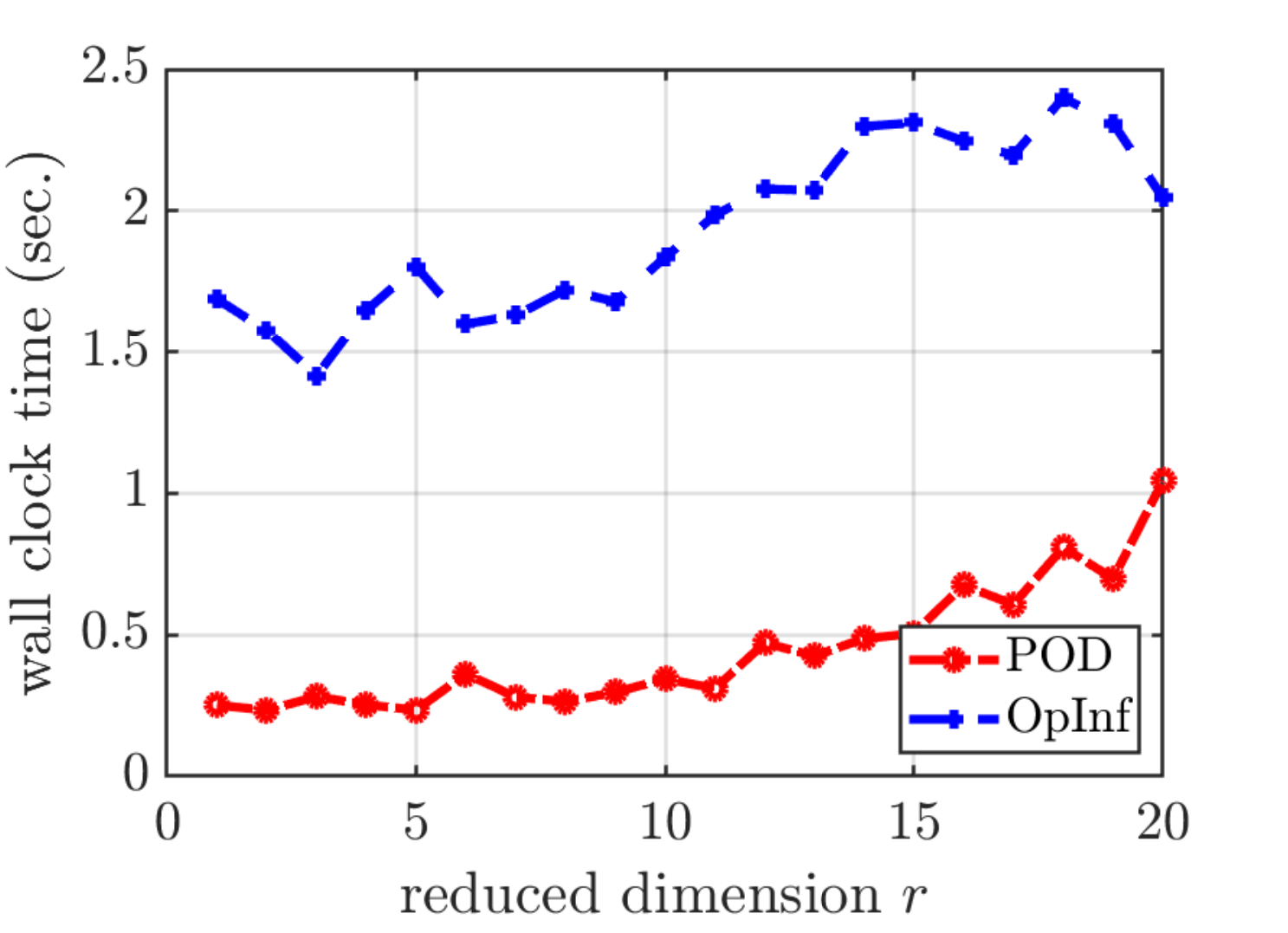}
	\caption{Elapsed total computational time of the ROMs. \label{fig:cpu}}	
\end{figure}

Next, we show the total computational times for the POD-G and OpInf algorithms in Figure \ref{fig:cpu}, where we exclude the computation time of obtaining POD basis since it is common for the POD-G and OpInf algorithms. The computational time consists of the time for constructing the reduced tensors \eqref{redquad} and online computation time for the POD-G algorithm, whereas in the OpInf algorithm it consists of the elapsed time of constructing the re-projected states and time derivatives, solving the least-squares problem \eqref{eq:lsqr}, and online computation time. The computational cost of the OpInf is larger than the POD-G due to re-projection. The total computational times in Figure \ref{fig:cpu} are  $1.044$s for the POD-G, and $2.045$s for  OpInf, both of which show a speedup of order $10^2$ over the FOM with $585.129$s computing time.

\begin{figure}[H]
	\centering
	\subfloat{\includegraphics[width=0.33\columnwidth]{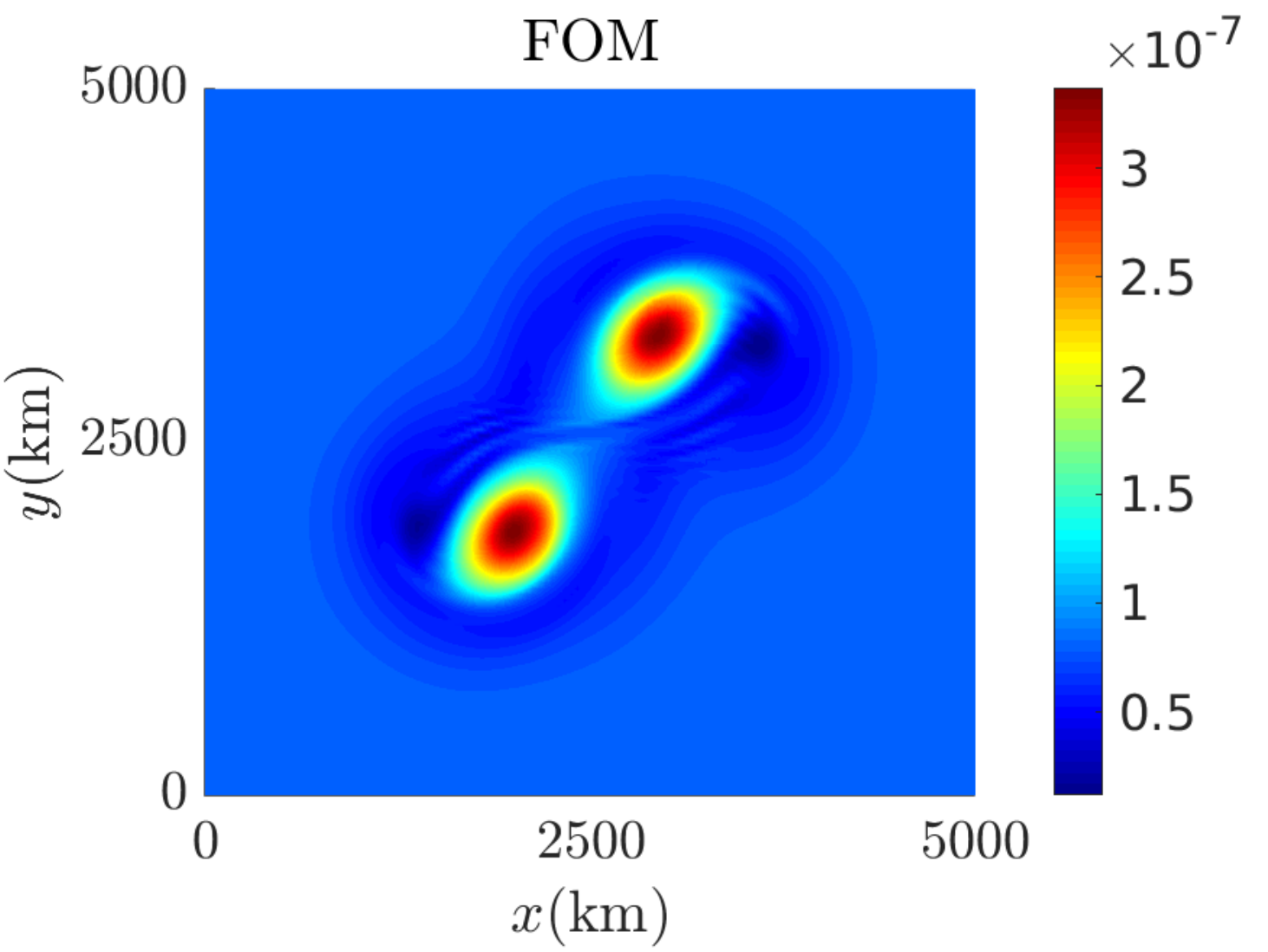}}
	\subfloat{\includegraphics[width=0.33\columnwidth]{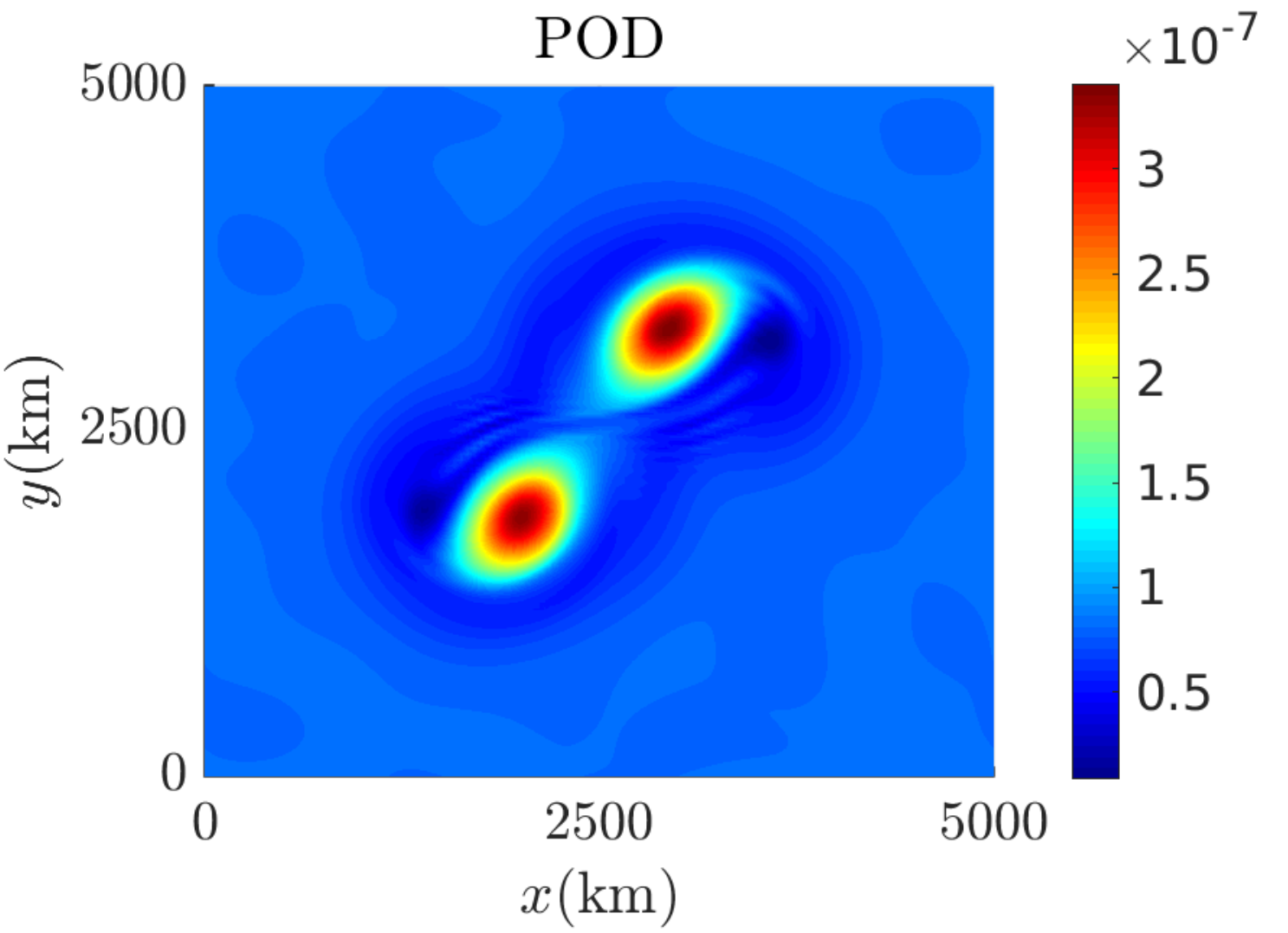}}
	\subfloat{\includegraphics[width=0.33\columnwidth]{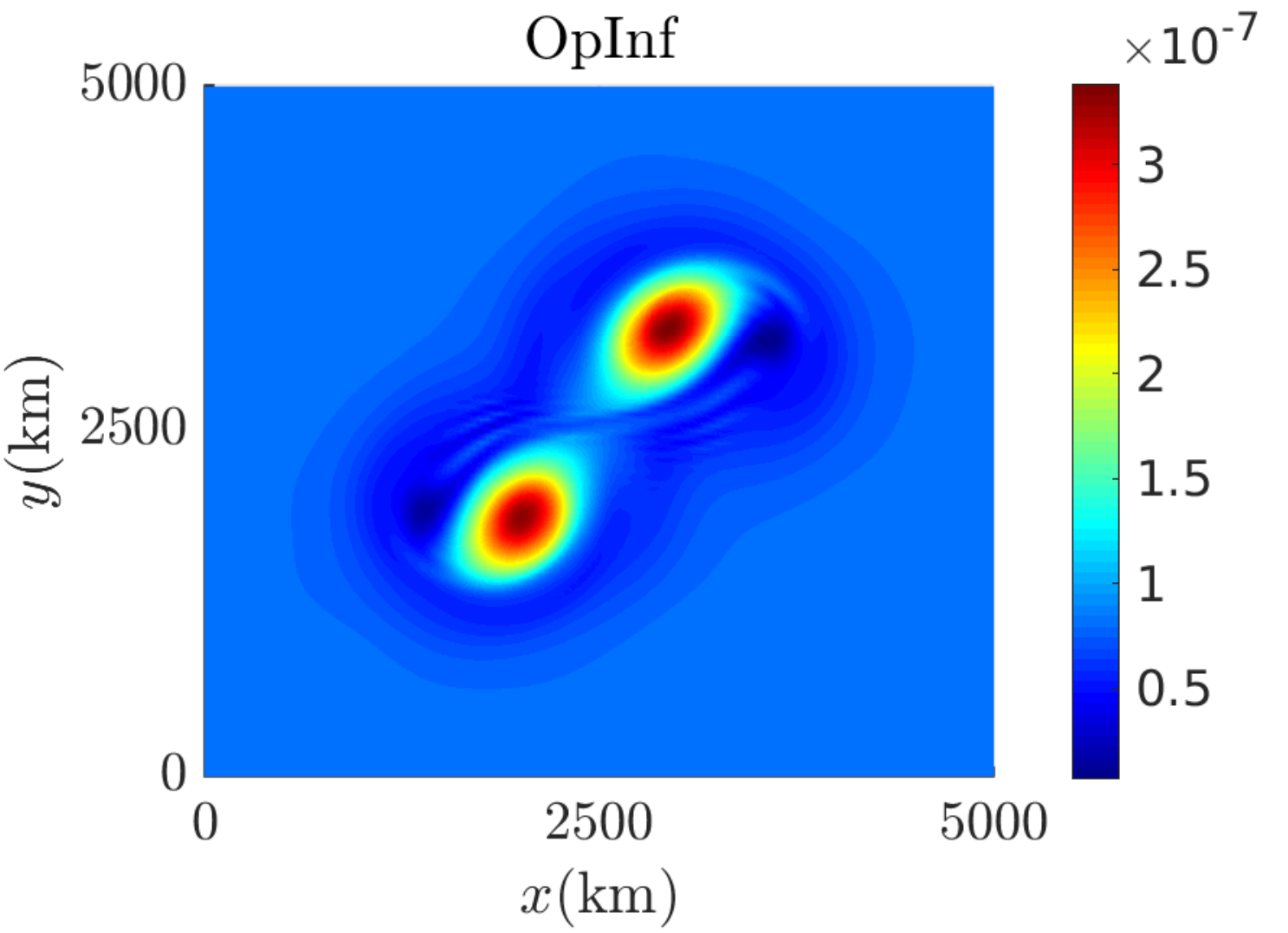}}
	
   \subfloat{\includegraphics[width=0.33\columnwidth]{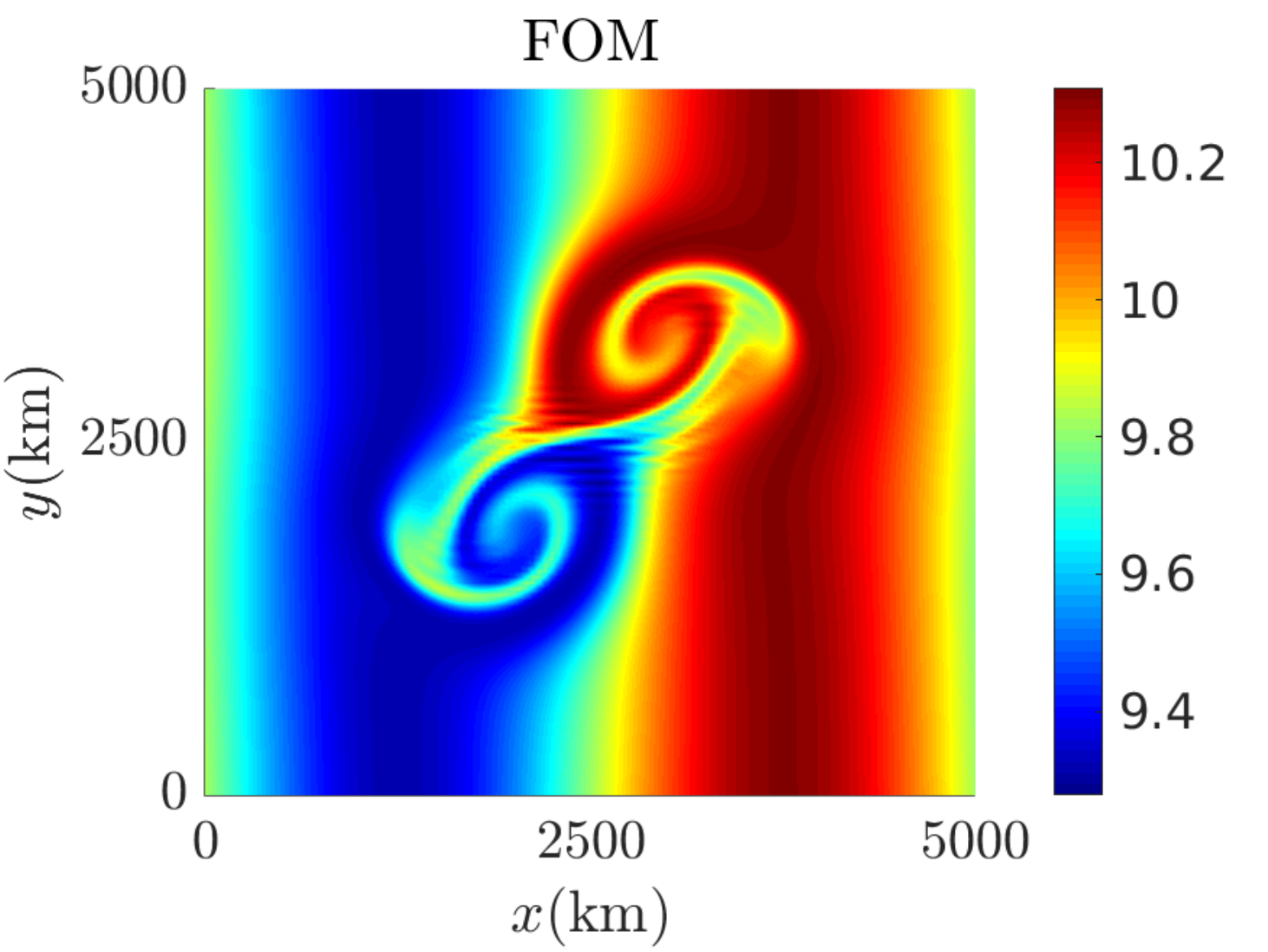}}
   \subfloat{\includegraphics[width=0.33\columnwidth]{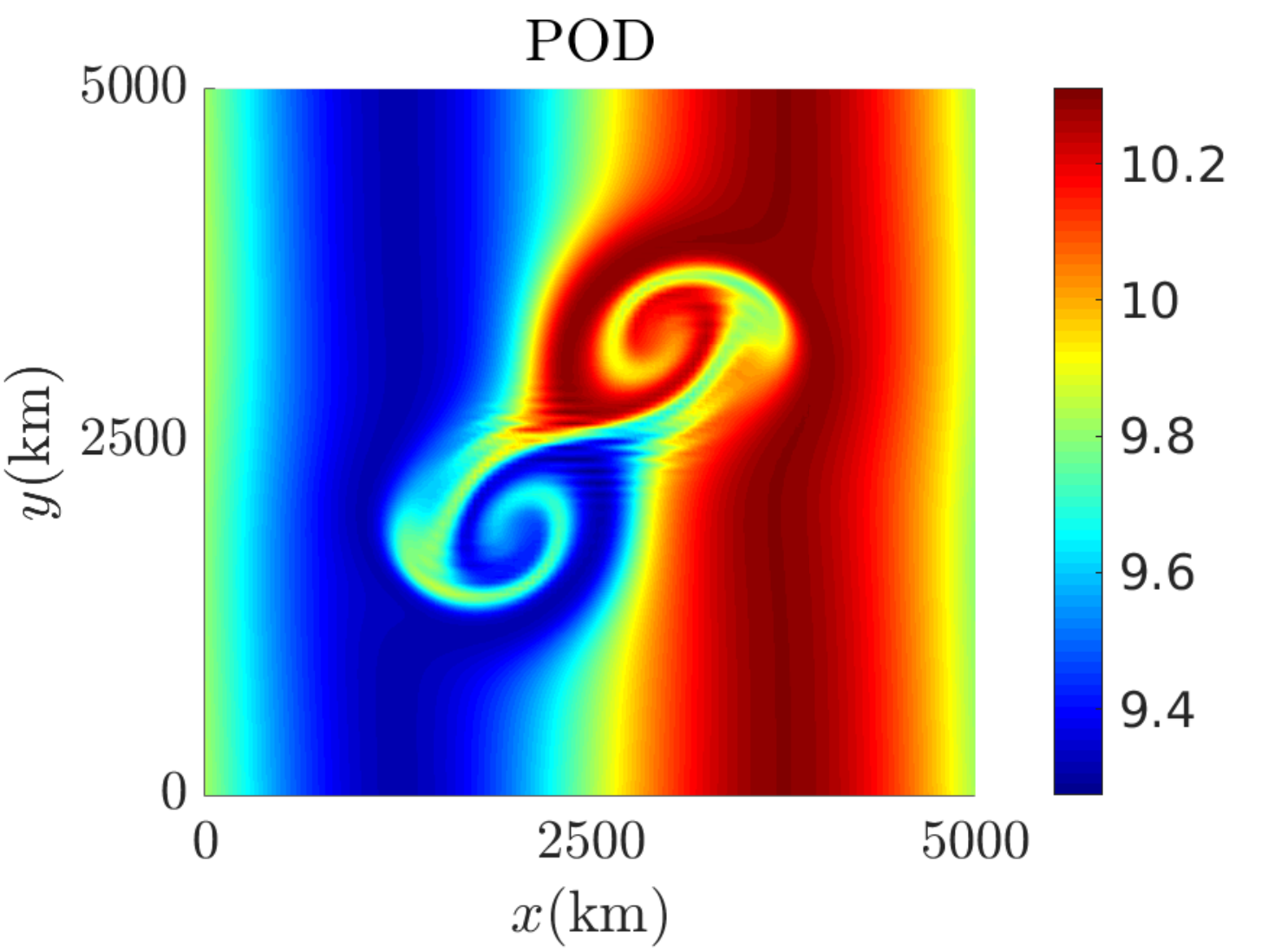}}
   \subfloat{\includegraphics[width=0.33\columnwidth]{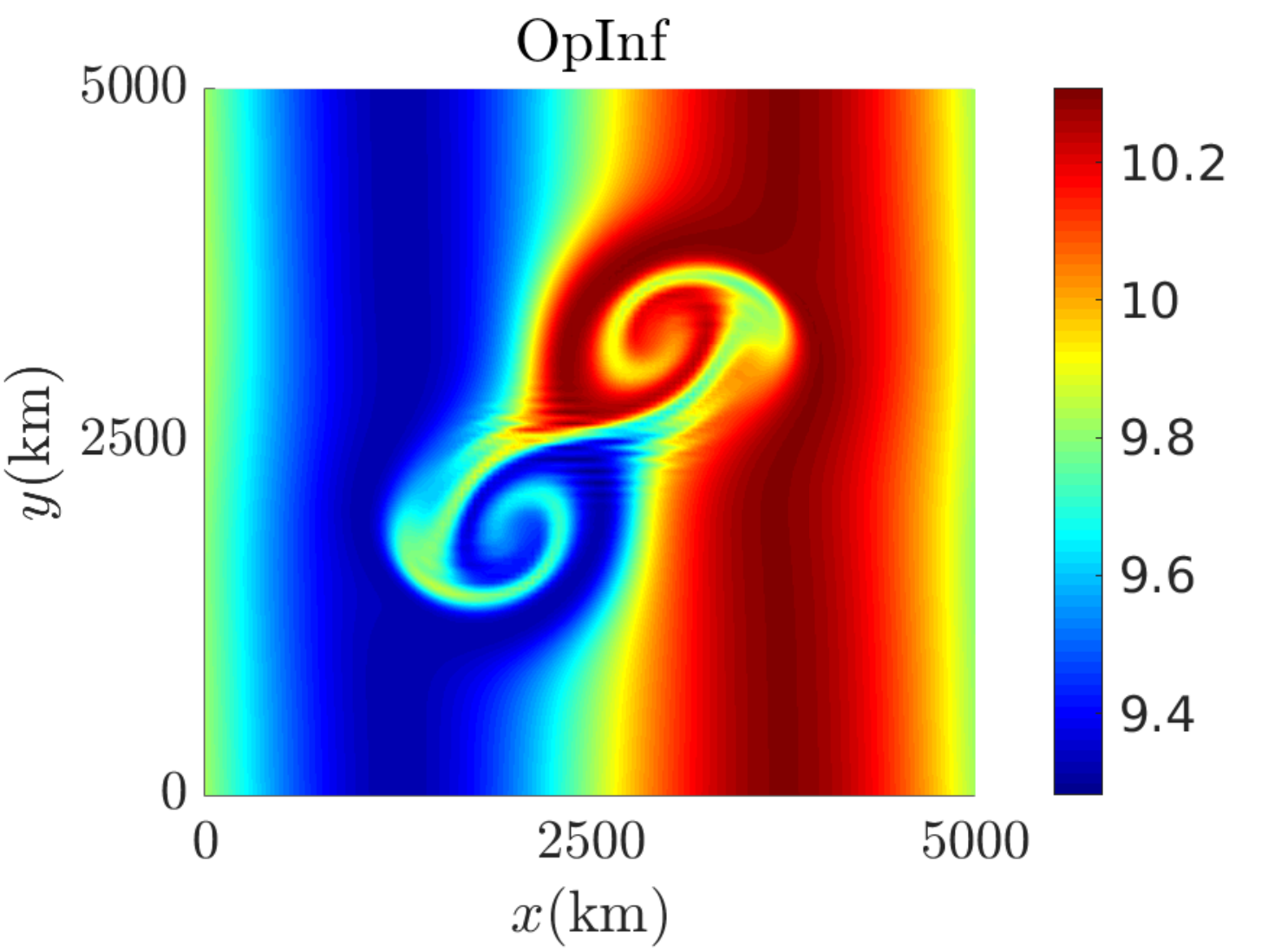}}
	\caption{(Top) Vorticity and (bottom) bouyancy of the FOM and ROMs at the final time. \label{fig:vort}}
\end{figure}

\begin{figure}[H]
	\centering
		\subfloat{\includegraphics[width=0.4\columnwidth]{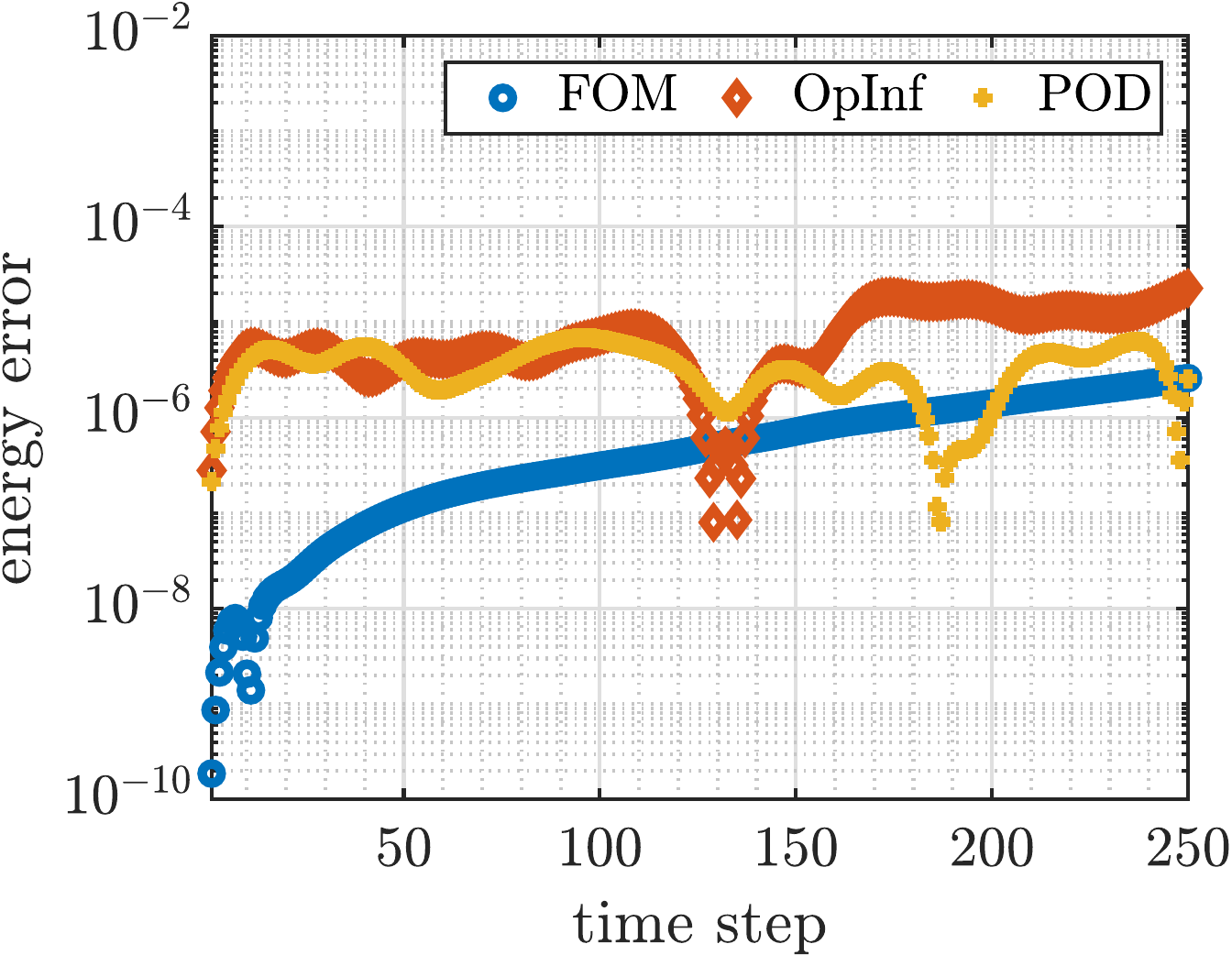}}
		\subfloat{\includegraphics[width=0.4\columnwidth]{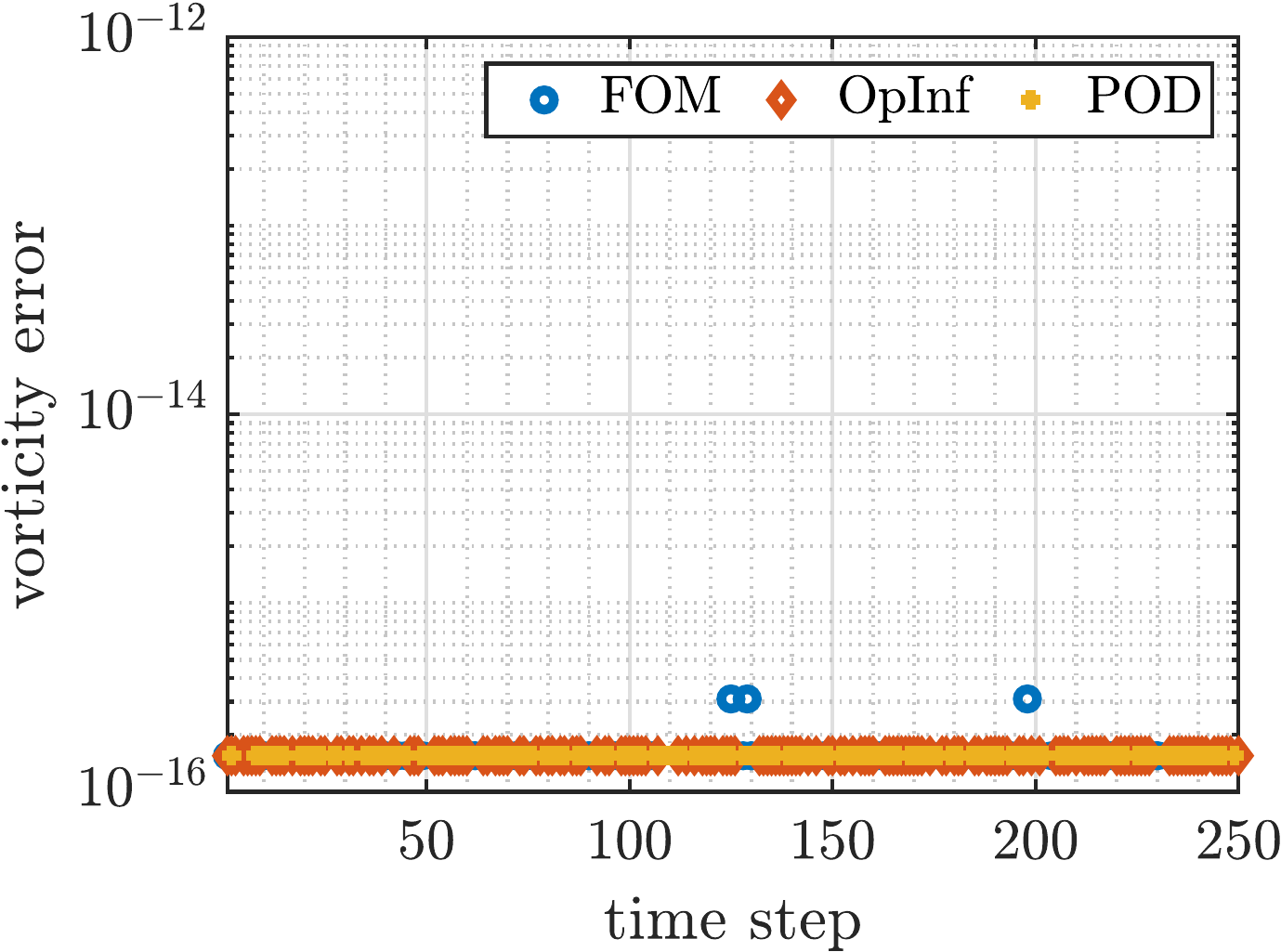}}
		
		\subfloat{\includegraphics[width=0.4\columnwidth]{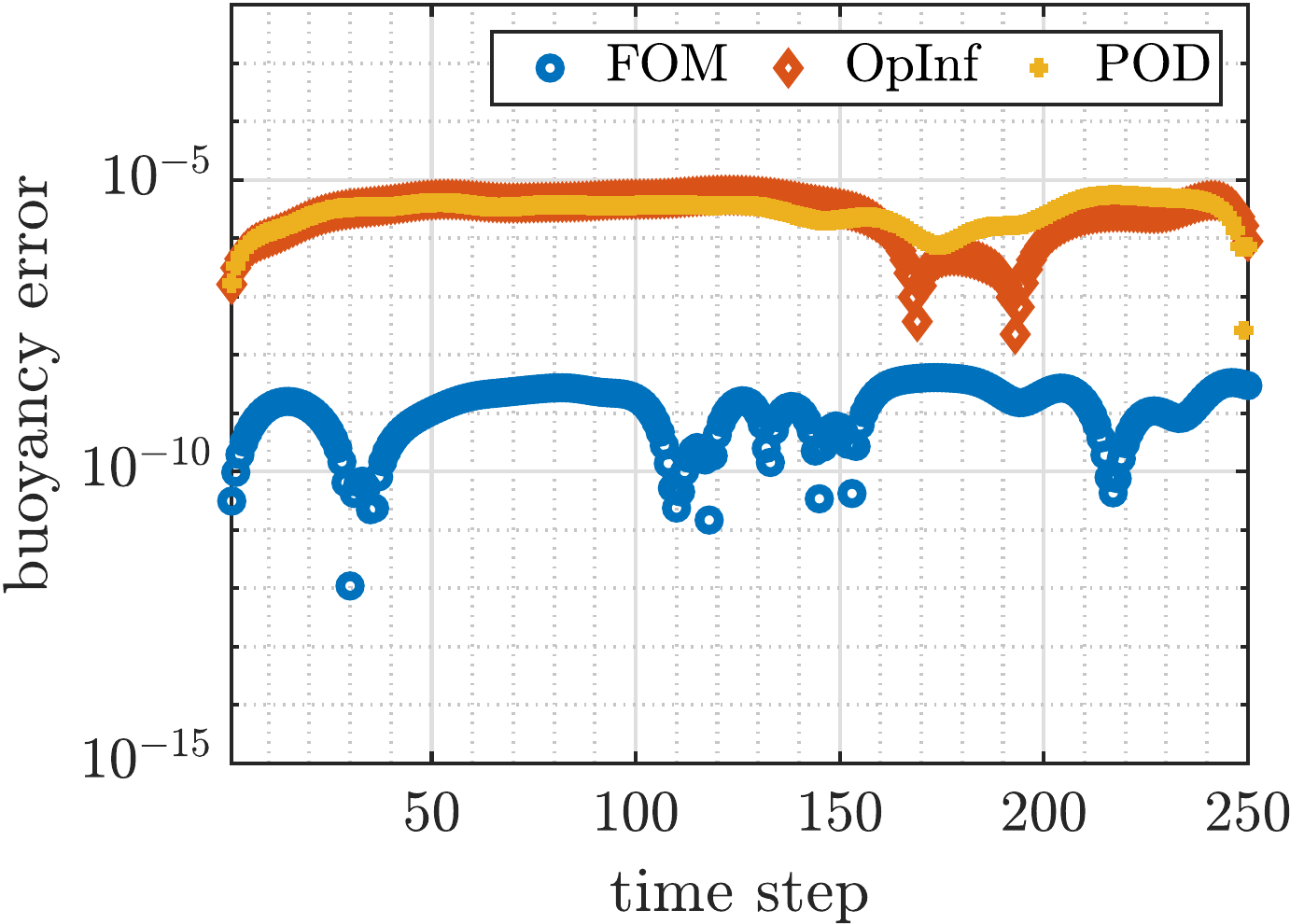}}
		\subfloat{\includegraphics[width=0.4\columnwidth]{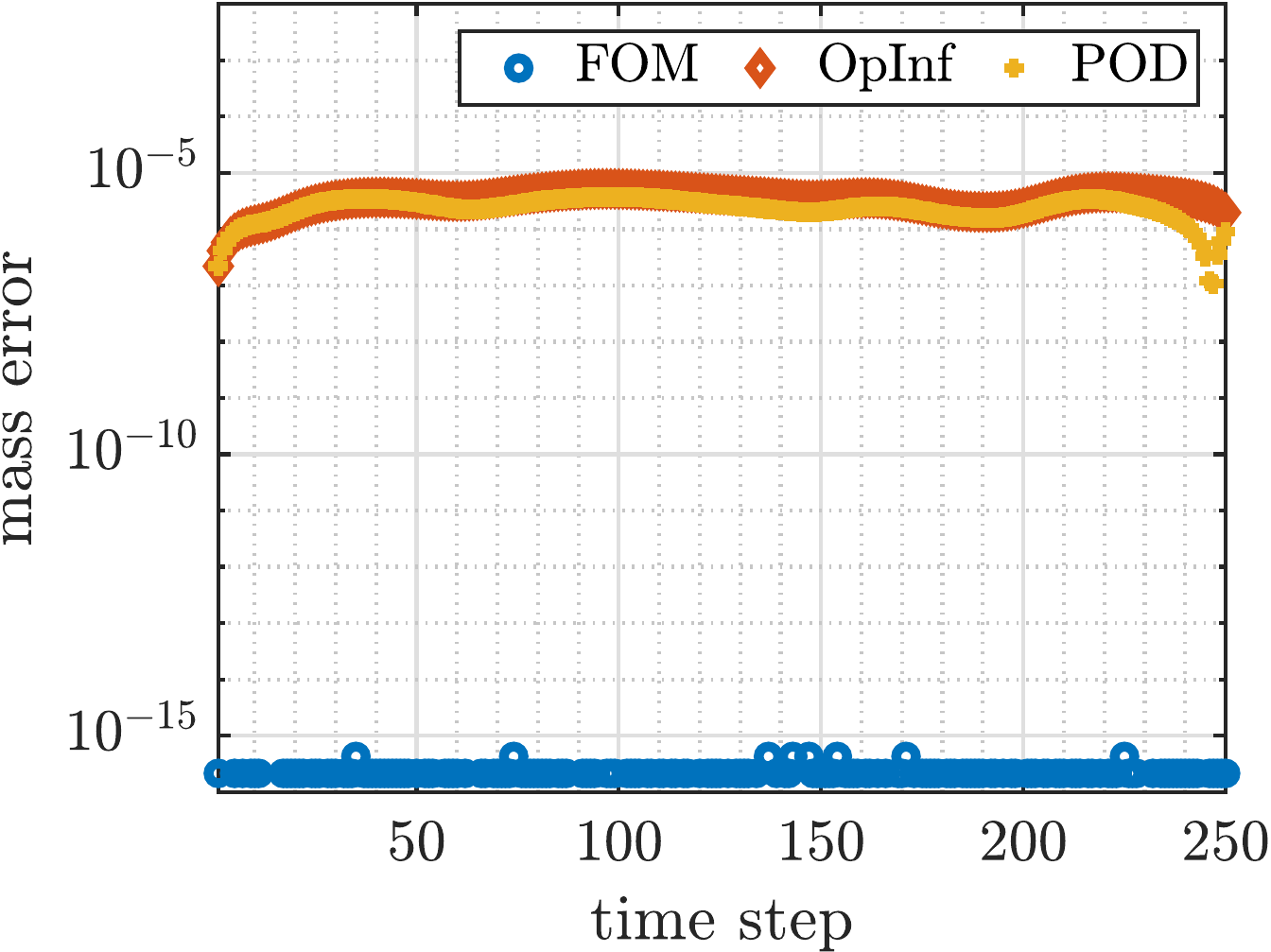}}
		\caption{Relative errors in the  conserved quantities. \label{conserved}}
\end{figure}

Both ROMs reproduce the full solutions accurately at the final time in Figure \ref{fig:vort} for the reduced dimension $r=20$. In Figure \ref{conserved}, the relative errors in the  Hamiltonian (energy)  $|H^k -H^0|/|H^0|$, mass $|M^k - M^0|/|M^0|$, buoyancy $|B^k - B^0|/|B^0|$, and the total potential vorticity  $|Q^k -Q^0|/|Q^0|$ are plotted over the time steps. The total mass and the total potential vorticity are preserved up to machine precision. The energy and buoyancy errors of both ROMs show bounded oscillations over time without any drift, i.e., they are preserved well in long-term. Preservation of conserved quantities by the ROMs is shown in Table~\ref{tbl2} in terms of the relative errors  \eqref{abserrorcon}.

\begin{table}[H]
	\centering
			\caption{Time-averaged relative errors \eqref{abserrorcon} of the conserved
				quantities. \label{tbl2}}
	\begin{tabular}{l|l|l|l|l}
		& Energy    & Total Vorticity & Mass      & Buoyancy  \\ \hline
	FOM   & 7.484e-07 & 6.392e-17 & 1.545e-16 & 1.567e-09 \\
	POD   & 3.489e-06 & 1.024e-16 & 2.489e-06 & 3.053e-06 \\
	OpInf & 8.114e-06 & 1.018e-16 & 3.440e-06 & 3.050e-06
	\end{tabular}
\end{table}

Prediction of  the future dynamics of complex systems are investigated  with the intrusive and  non-intrusive ROMs  for the single-injector combustor \cite{Willcox20}, the SWE \cite{Navon20}, and  the quasi-geostrophic equations \cite{Iliescu20}. 
Here, we demonstrate the predictive capabilities of  both ROM approaches for the RTSWE.
For the state vector  ${\bm w}= ({\bm u}, {\bm v},{\bm h}, {\bm s} )$, the relative FOM-ROM error at the time step $k$ is defined as
\begin{align}\label{relerrts}
\frac{\|{\bm  w}^k-{\bm  w}_r^k\|_{L^2}}{\|{\bm  w}^k\|_{L^2}}.
\end{align}

\begin{figure}[H]
\centering
\includegraphics[width=0.45\columnwidth]{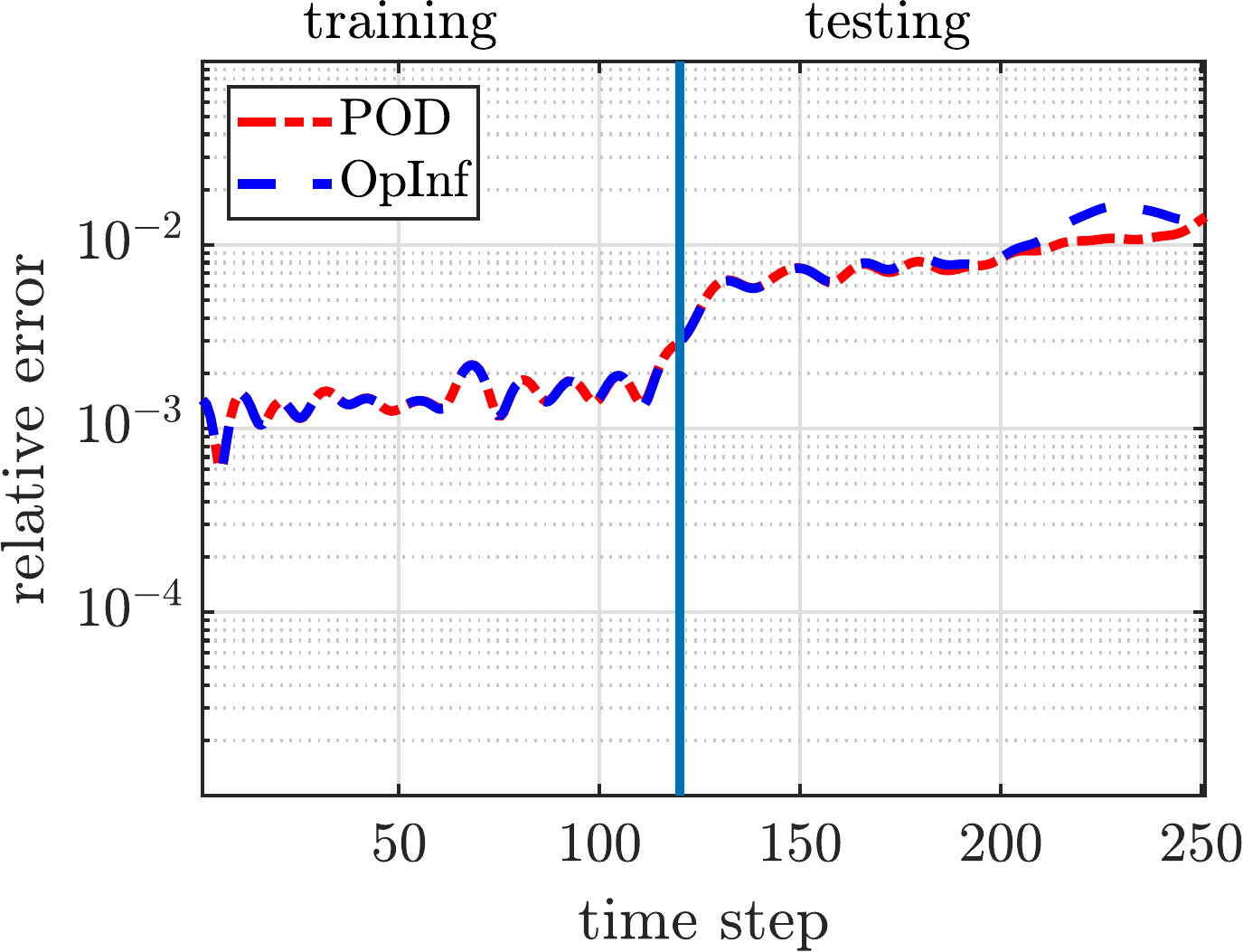}
\includegraphics[width=0.45\columnwidth]{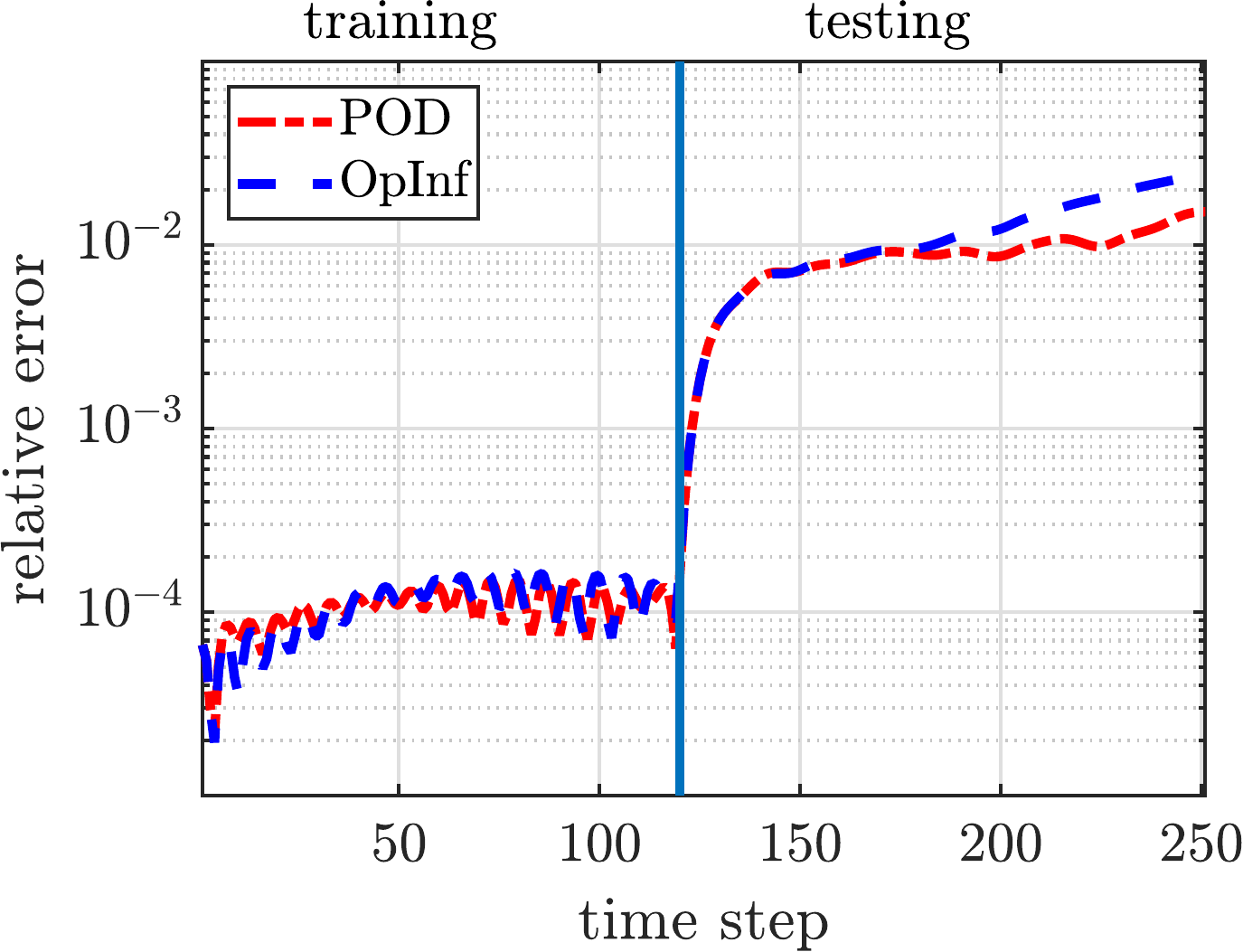}
	\caption{Prediction performance of the ROMs trained up to $ K=120 $ with the reduced dimension (left) $ r=10 $ and (right) $ r=20 $. \label{fig:predpods}}
\end{figure}

\begin{figure}[H]
\centering
\subfloat{\includegraphics[width=0.45\columnwidth]{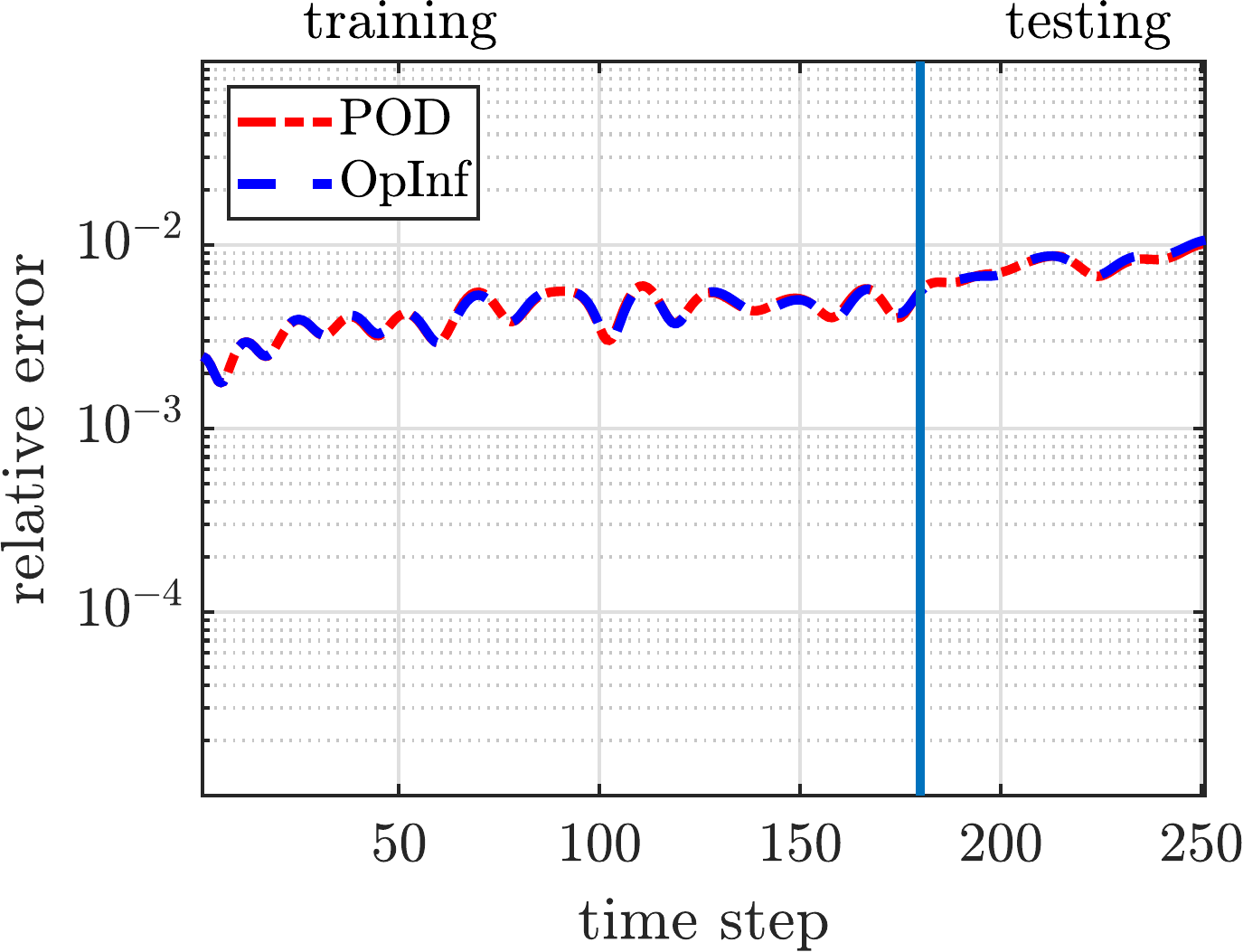}}
\subfloat{\includegraphics[width=0.45\columnwidth]{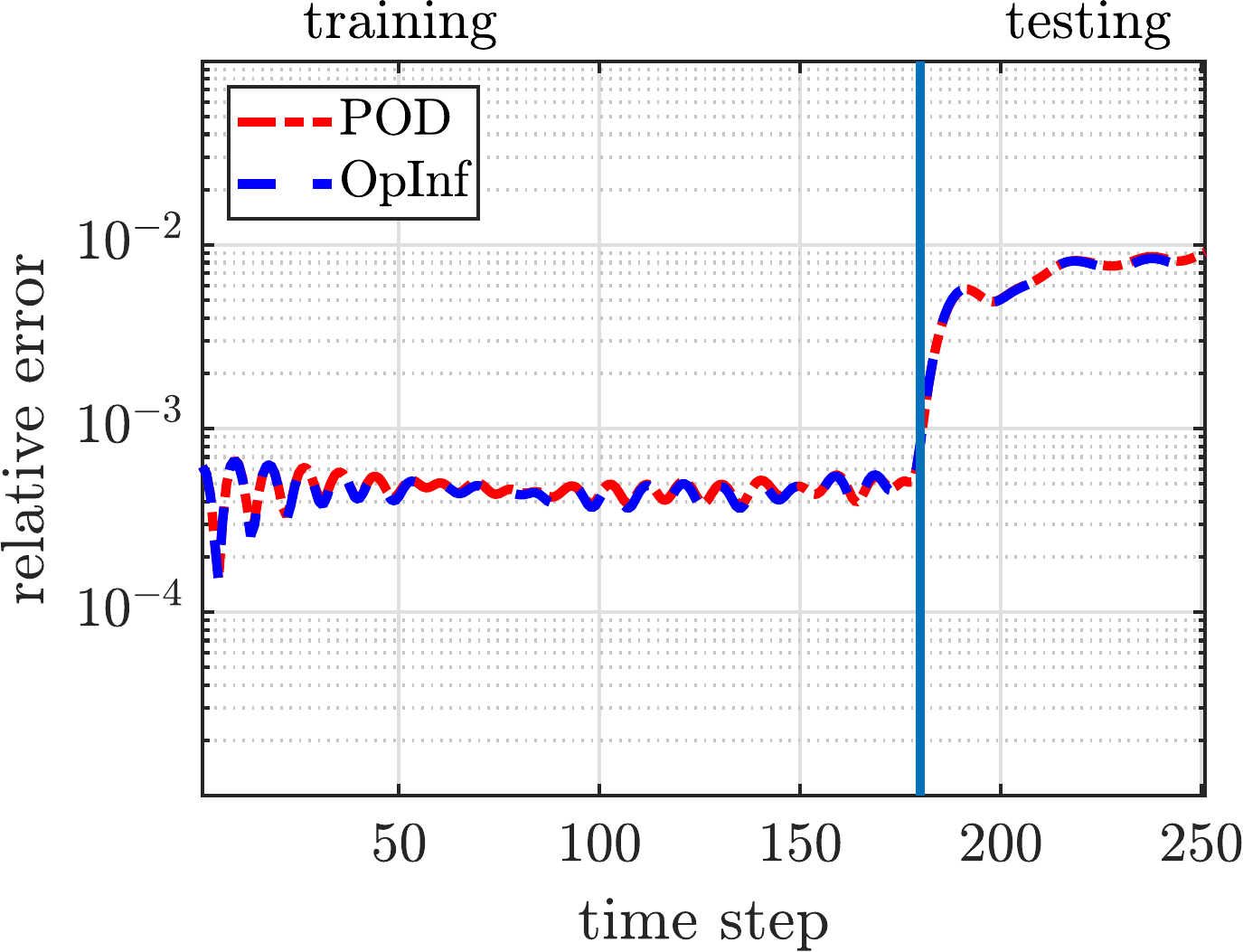}}
	\caption{Prediction performance of the ROMs trained up to $ K=120 $ with the reduced dimension (left) $ r=10 $ and (right) $ r=20 $. \label{fig:predpod10}}
\end{figure}

In Figures \ref{fig:predpods}-\ref{fig:predpod10}, the relative errors  \eqref{relerrts} for the state variables are plotted, where the vertical blue lines separate the training  and the prediction periods.
The training period $[t_0,t_{k-1}]$ is taken of the full time domain $[t_0,t_f]$ for prediction of the  dynamics of the buoyancy and the vorticity over the whole time domain $[t_k,t_f]$. We study two cases; shorter and longer   training  regimes than the prediction regimes in Figures  \ref{fig:predpods}  and \ref{fig:predpod10}.
When the amount of training data gets larger than the prediction period, the ROM predictions improve in Figure~\ref{fig:predpod10}. Increasing the number of the modes does not affect the ROM predictions significantly in Figures \ref{fig:predpods}-\ref{fig:predpod10}. In each case, the ROMs are able to accurately re-predict the training data and capture much of the overall system behavior in the prediction phase.
Both POD-G and OpInf are able to accurately re-predict the training data and capture much of the overall behavior of the FOMs in  the prediction period. An overview of the performance of ROMs in terms  of the averaged relative errors \eqref{relerrorsol} is given in Table  \eqref{table:pred}.

\begin{table}[H]
	\caption{Average relative errors \eqref{relerrorsol} for training and prediction sets. \label{table:pred}}
	\centering
	\begin{tabular}{l|l|l|l|l|l}
		\multicolumn{2}{c|}{}            & \multicolumn{2}{c|}{$ r=10 $} & \multicolumn{2}{c}{$ r=20 $} \\ \hline
		\multicolumn{2}{c|}{}            & POD         & OpInf       & POD         & OpInf       \\ \hline
		\multicolumn{1}{c|}{\multirow{2}{*}{$ K=120 $}} & training   & 1.529e-03   & 1.523e-03    & 1.060e-04   & 1.106e-04    \\ \cline{2-6}
		\multicolumn{1}{c|}{}                       & prediction & 8.299e-03     & 9.487e-03    & 8.769e-03    & 1.193e-02   \\ \hline
		\multirow{2}{*}{$ K=180 $}                       & training   & 4.250e-03   & 4.233e-03   & 4.737e-04   & 4.637e-04   \\ \cline{2-6}
		& prediction & 7.691e-03   & 7.817e-03   & 6.695e-03   & 6.587e-03    \\
	\end{tabular}
\end{table}

\subsection{The parametric case }

For the parametric case, we consider the Coriolis parameter $f(\mu) = 2\Omega\sin(\mu)$ varying with the latitude $\mu$  in the parameter domain ${\mathcal D }: = 4\pi/18 \le  \mu \le 8\pi/18$, i.e., between the $40^0$th to $80^0 $th latitude.  We take a larger time horizon, with the number of time steps $ K=300 $, than the one in the non-parametric case in order to reflect the parametric variability of the reduced solutions and to improve the condition number of the data matrices.
The average relative errors in the training phase are computed  as
\begin{equation}\label{RelAvgError}
\frac{1}{M_{\mathrm train}}\sum_{i = 1}^{M_{\mathrm train}} \frac{\|\bm \Phi_r Z(\mu_i^{\mathrm train}) -
W(\mu_i^{\mathrm train})\|_F}{\|W(\mu_i^{\mathrm train})\|_F},
\end{equation}
where $ W(\mu_i^{\mathrm train}) $ is the FOM trajectory and ${\bm Z}(\mu_i^{\mathrm train})$ is the trajectory of either the intrusive reduced model (POD-G) or learned model (OpInf) from the re-projected trajectories.
Similarly, the accuracy of the parametric model for the test case is measured via the average relative error defined by
\begin{equation}\label{RelAvgErrortest}
\frac{1}{M^{\mathrm test}}\sum_{i = 1}^{M_{\mathrm test}} \frac{\|\bm \Phi_r Z(\mu_i^{\mathrm test}) -
	W(\mu_i^{\mathrm test})\|_F}{\|W(\mu_i ^{\mathrm test})\|_F}.
\end{equation}
In the training phase, the intrusive and non-intrusive ROMs  are  constructed with $M_{\mathrm train}= 6$ equidistant parameters $\mu_1^{\mathrm train},\dots,  \mu_{M_{\mathrm train}}^{\mathrm train} \in \mathcal{ D}$ in the parameter domain $\mathcal{D}$. Moreover, in order to reduce the condition number of data matrices, we take initial conditions that are randomly perturbed at the center of the vortices as $oy=0.1+\gamma$, where $\gamma$ is a uniformly distributed random value. They do not necessarily have a physical meaning,  they reflect the behavior of the FOM over many different initial conditions with the aim to provide more complete information about the FOM and reduce the condition number of the data matrices in the OpInf. The least-squares problem \eqref{eq:lsqr} is solved with \texttt{lsqminnorm} algorithm with the tolerance $ tol=1e-10 $, and by considering subset of the snapshots, i.e., every $2$nd snapshot.

\begin{figure}[htb!]
\centering
\subfloat{\includegraphics[width=0.45\columnwidth]{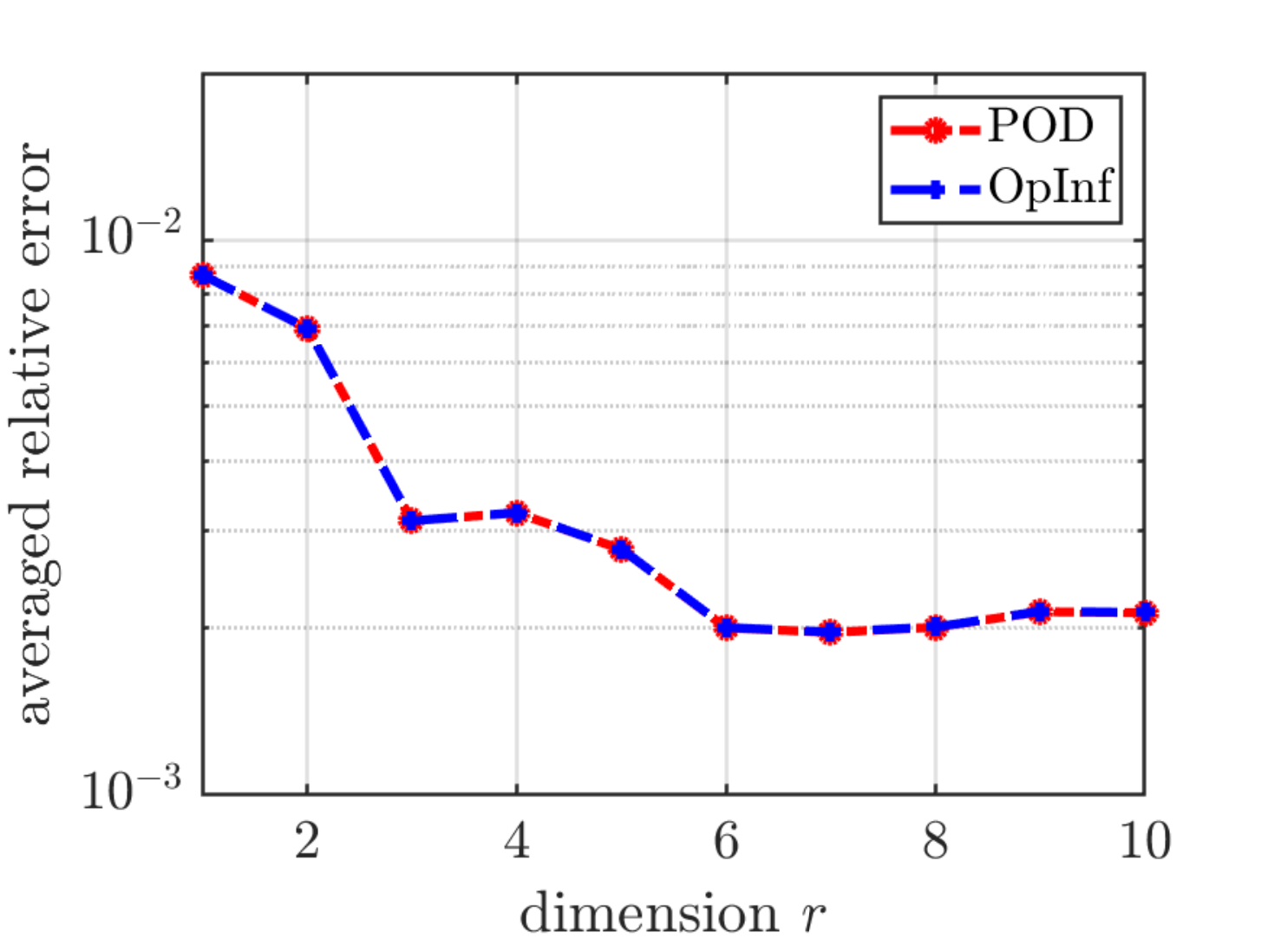}}
\subfloat{\includegraphics[width=0.45\columnwidth]{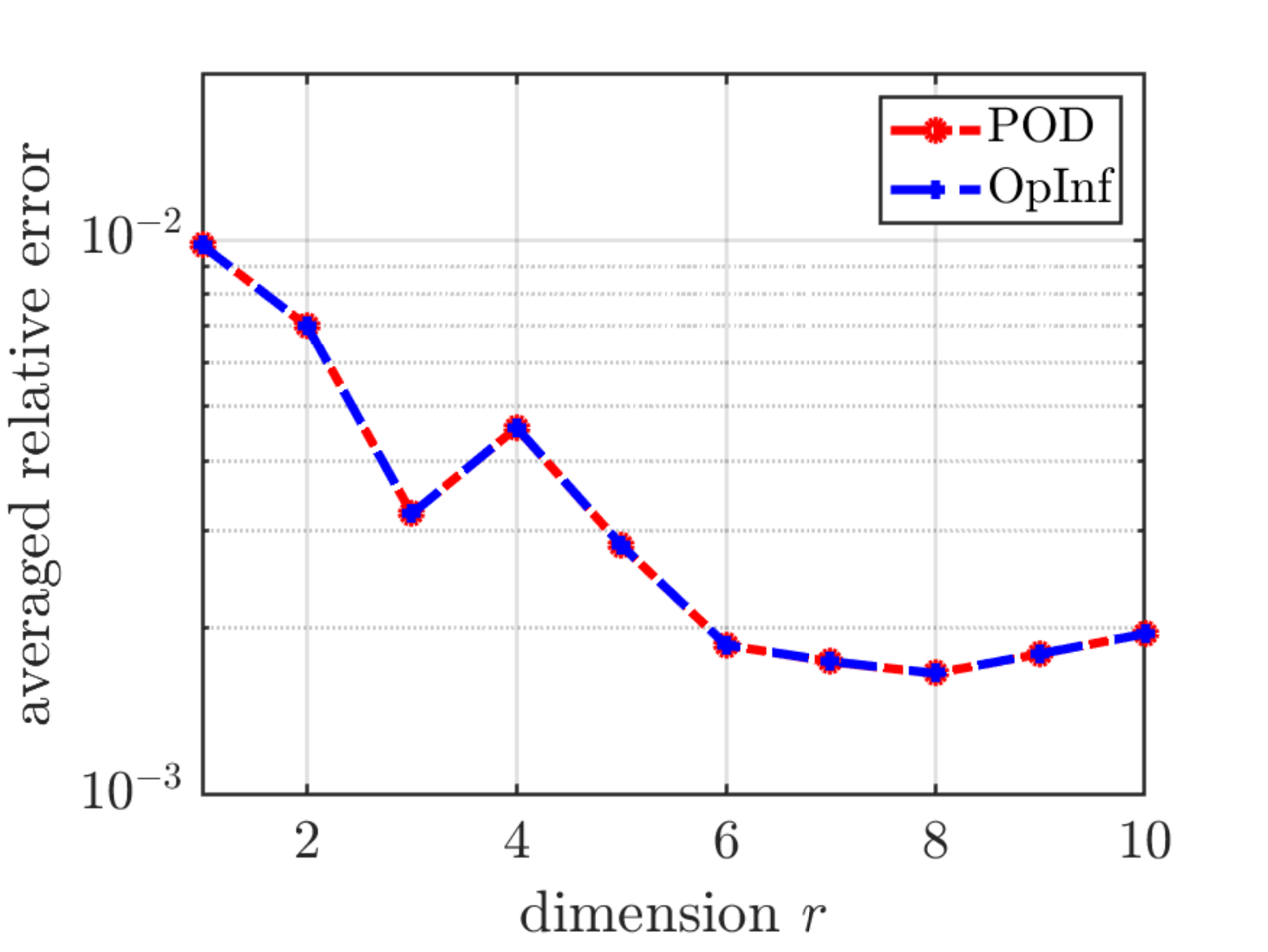}}
	\caption{(Left) The relative error \eqref{RelAvgError} in the training period, (right) the relative error \eqref{RelAvgErrortest} in the testing  period. \label{fig:par_rel_error}}
\end{figure}

\begin{figure}[htb!]
	\centering
	\subfloat{\includegraphics[width=0.33\columnwidth]{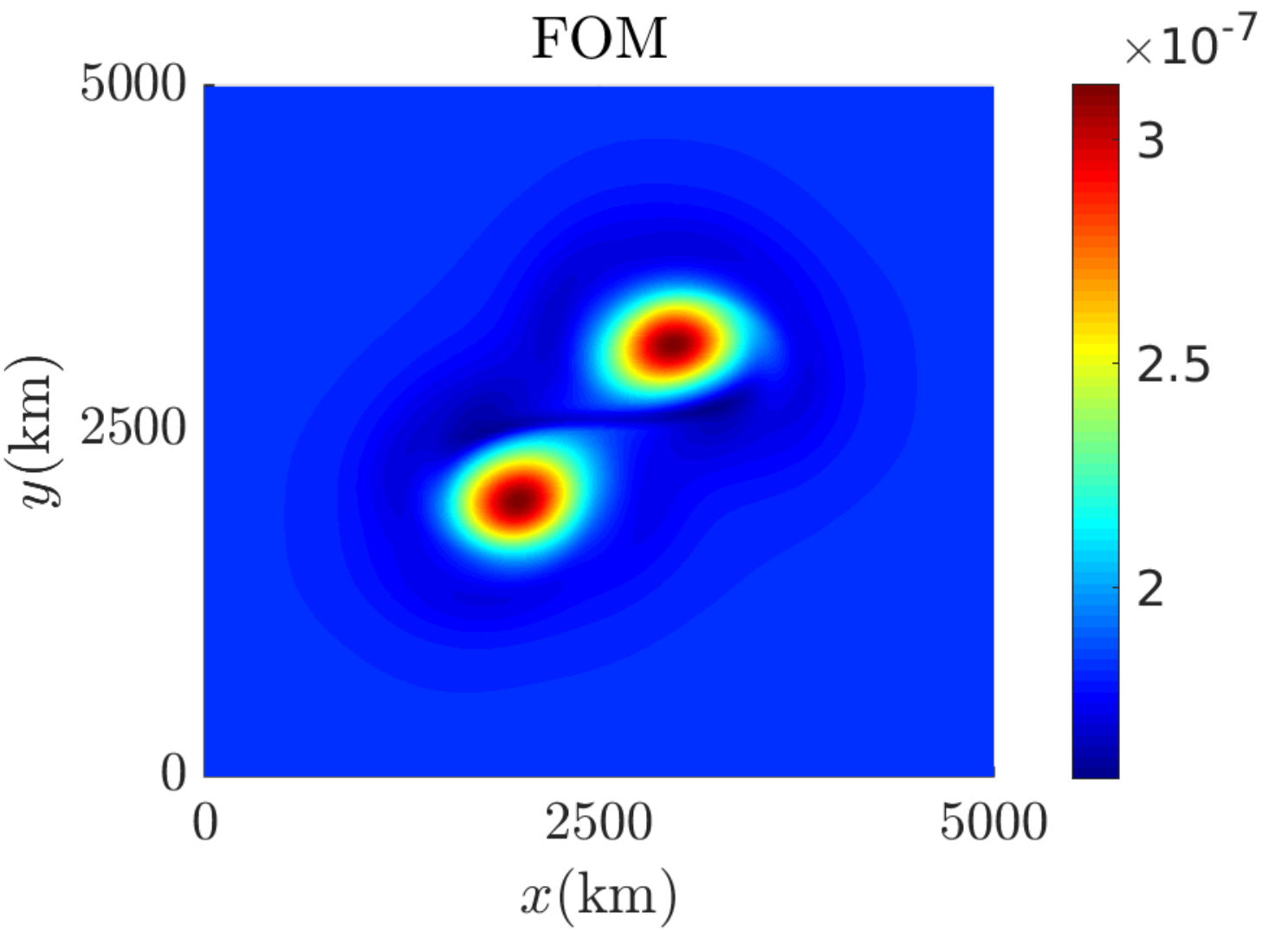}}	
	\subfloat{\includegraphics[width=0.33\columnwidth]{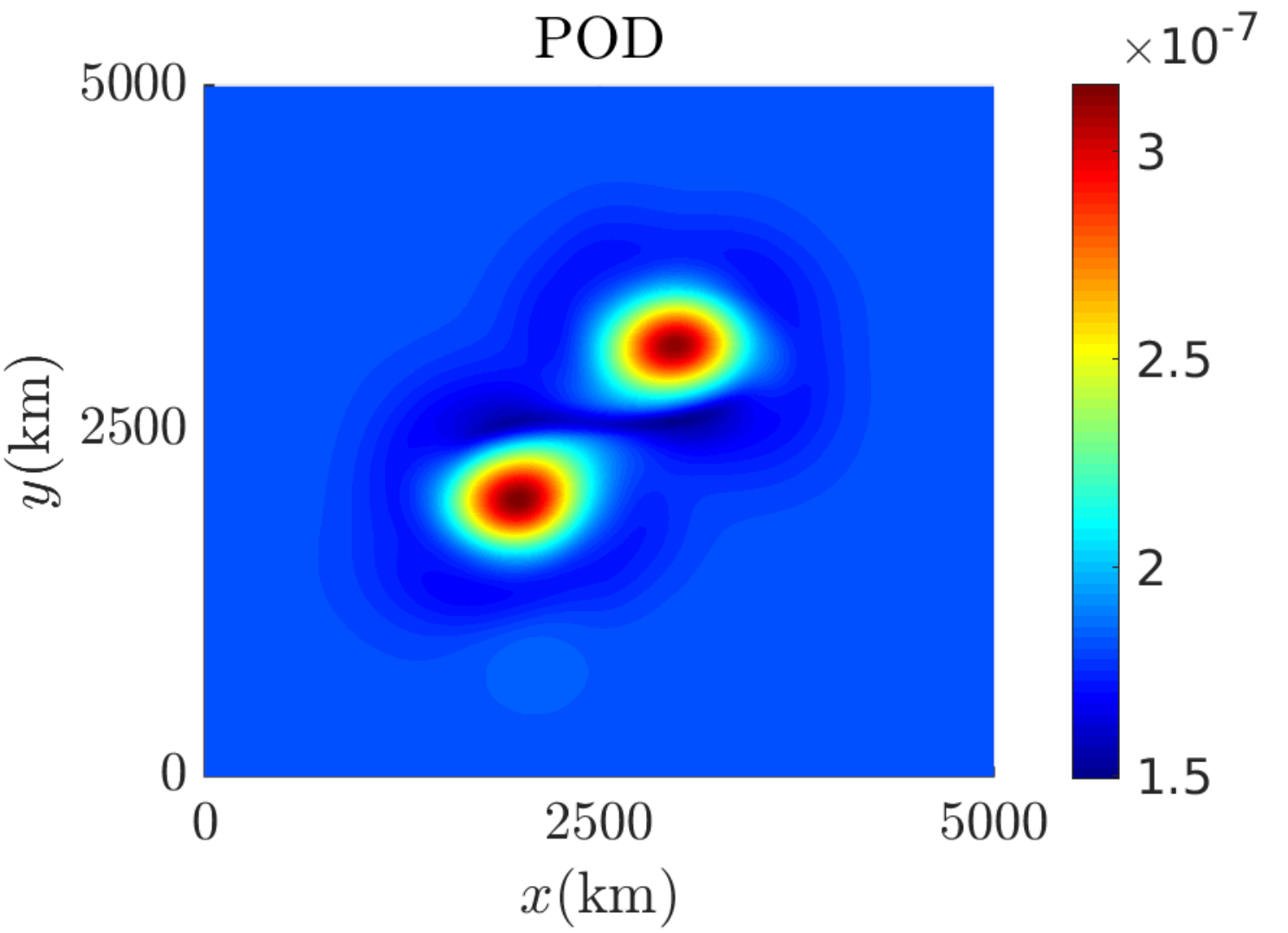}}
	\subfloat{\includegraphics[width=0.33\columnwidth]{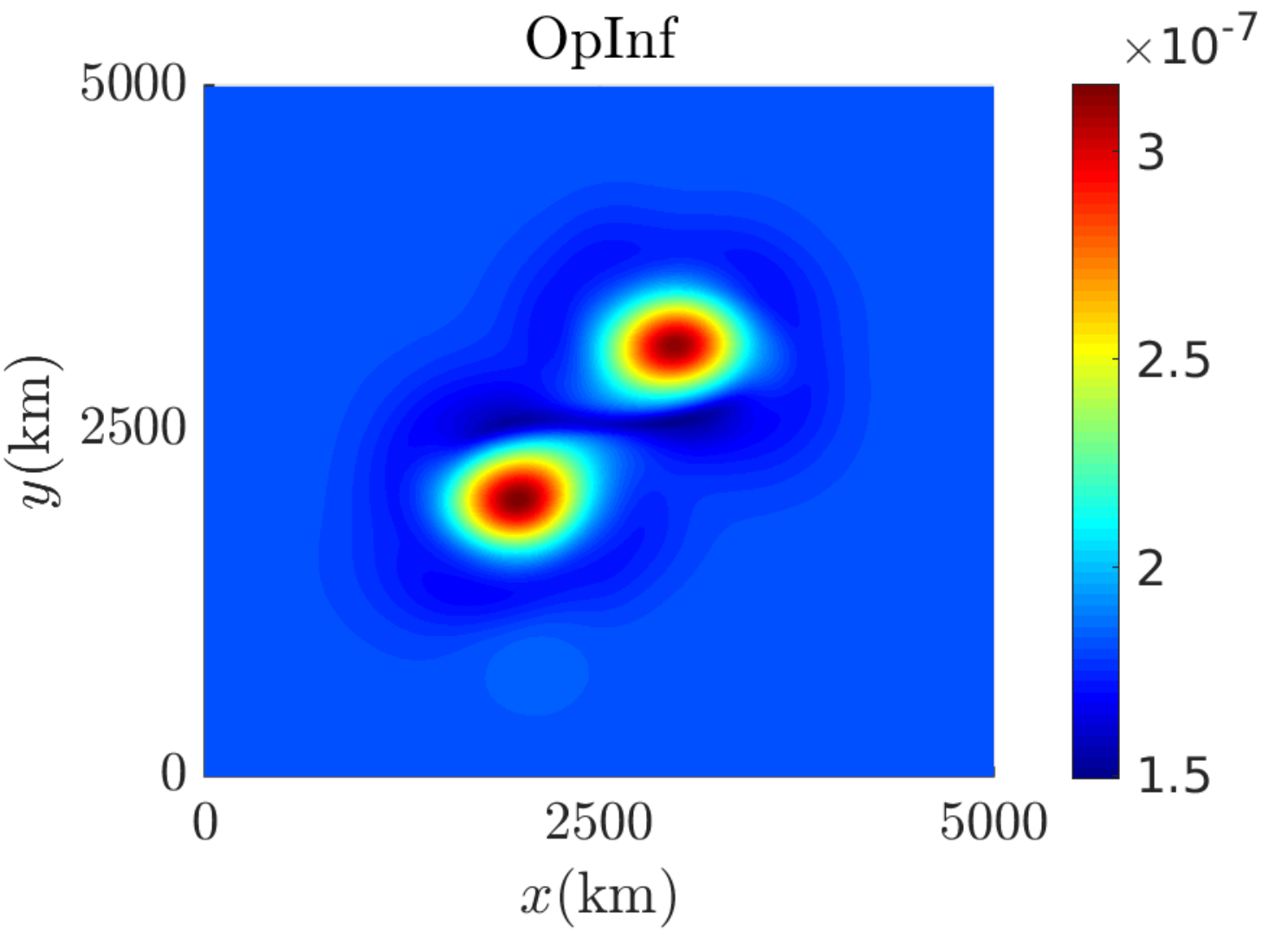}}
	
  \subfloat{\includegraphics[width=0.33\columnwidth]{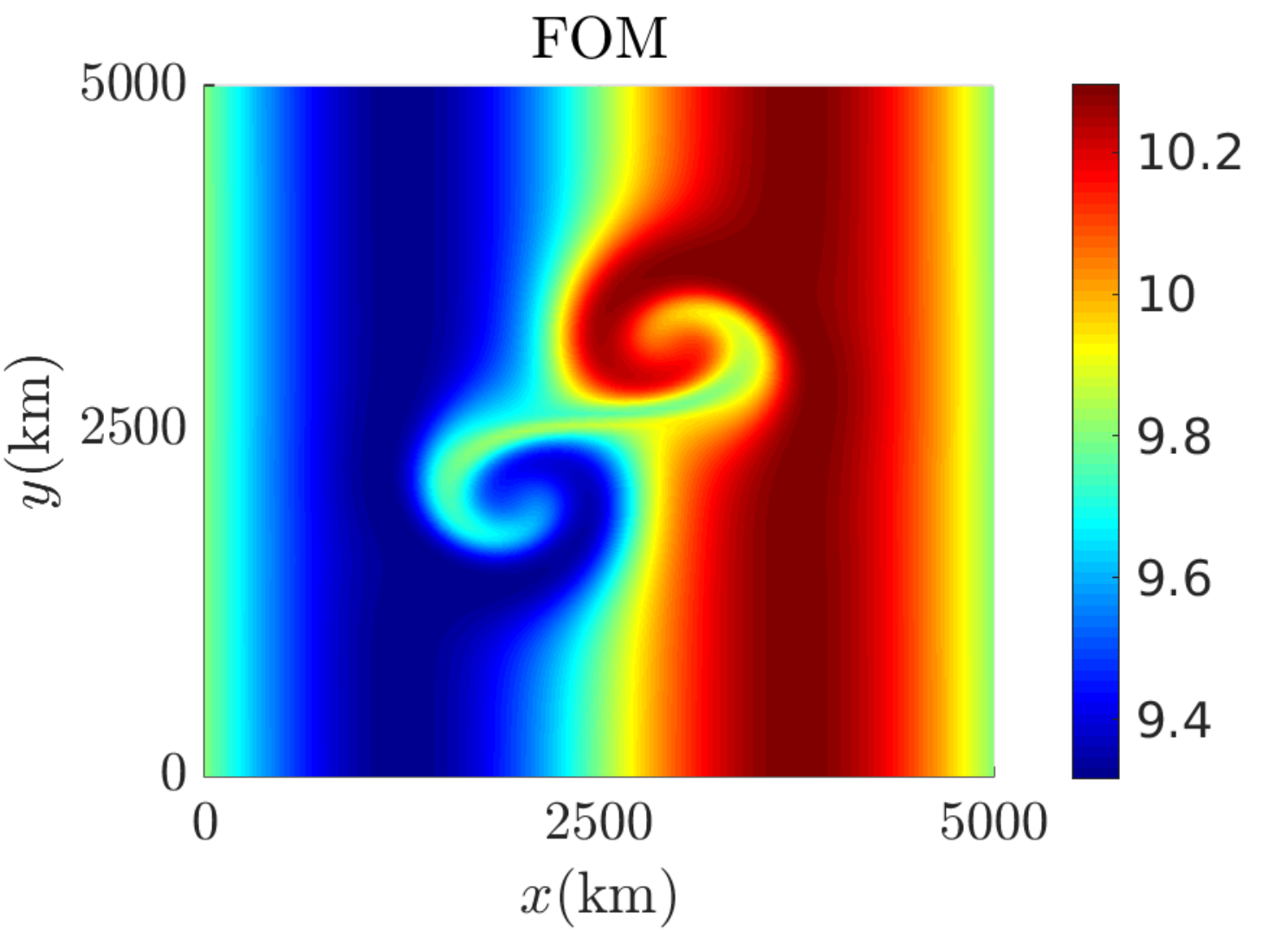}}	
  \subfloat{\includegraphics[width=0.33\columnwidth]{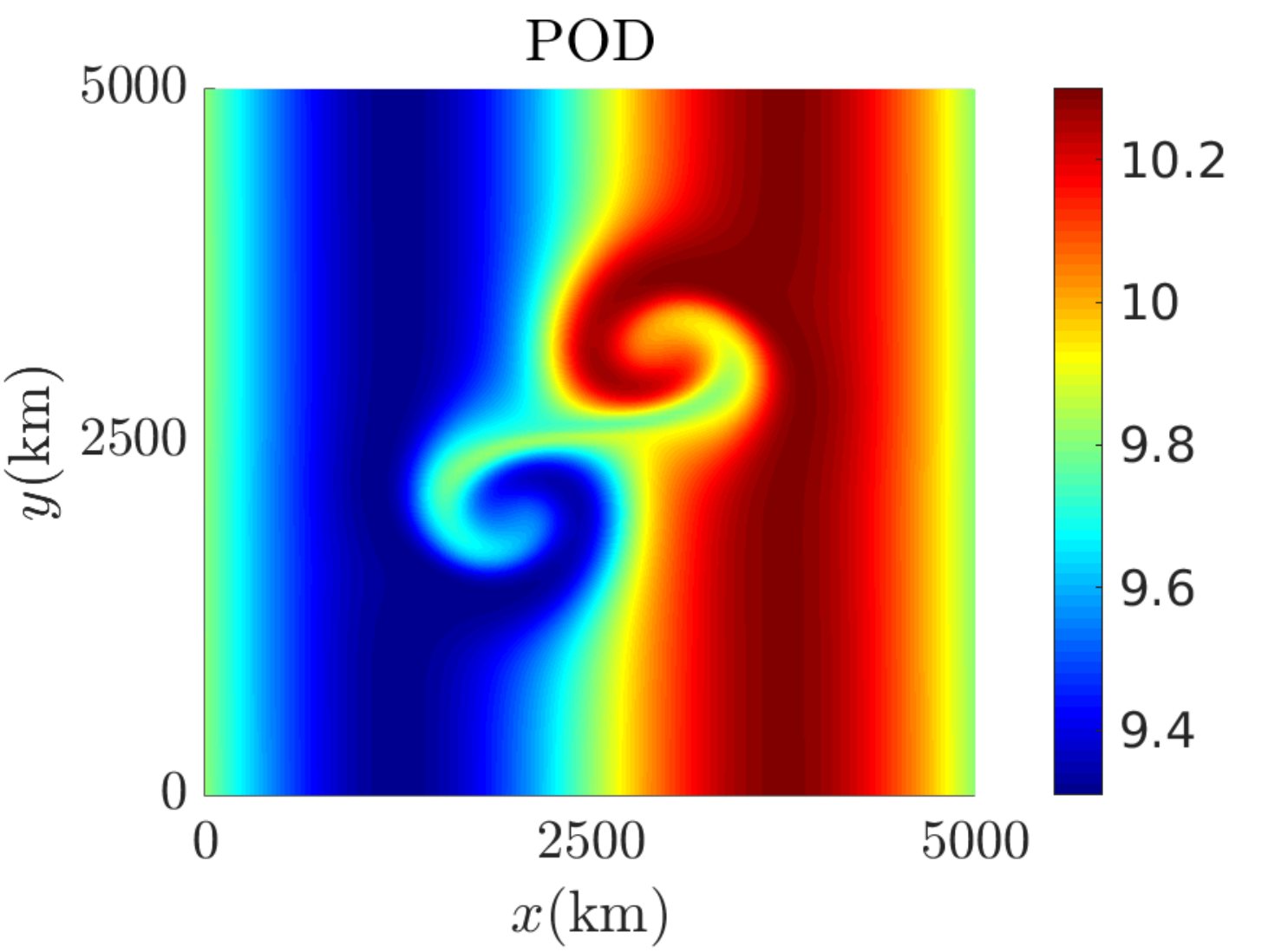}}
  \subfloat{\includegraphics[width=0.33\columnwidth]{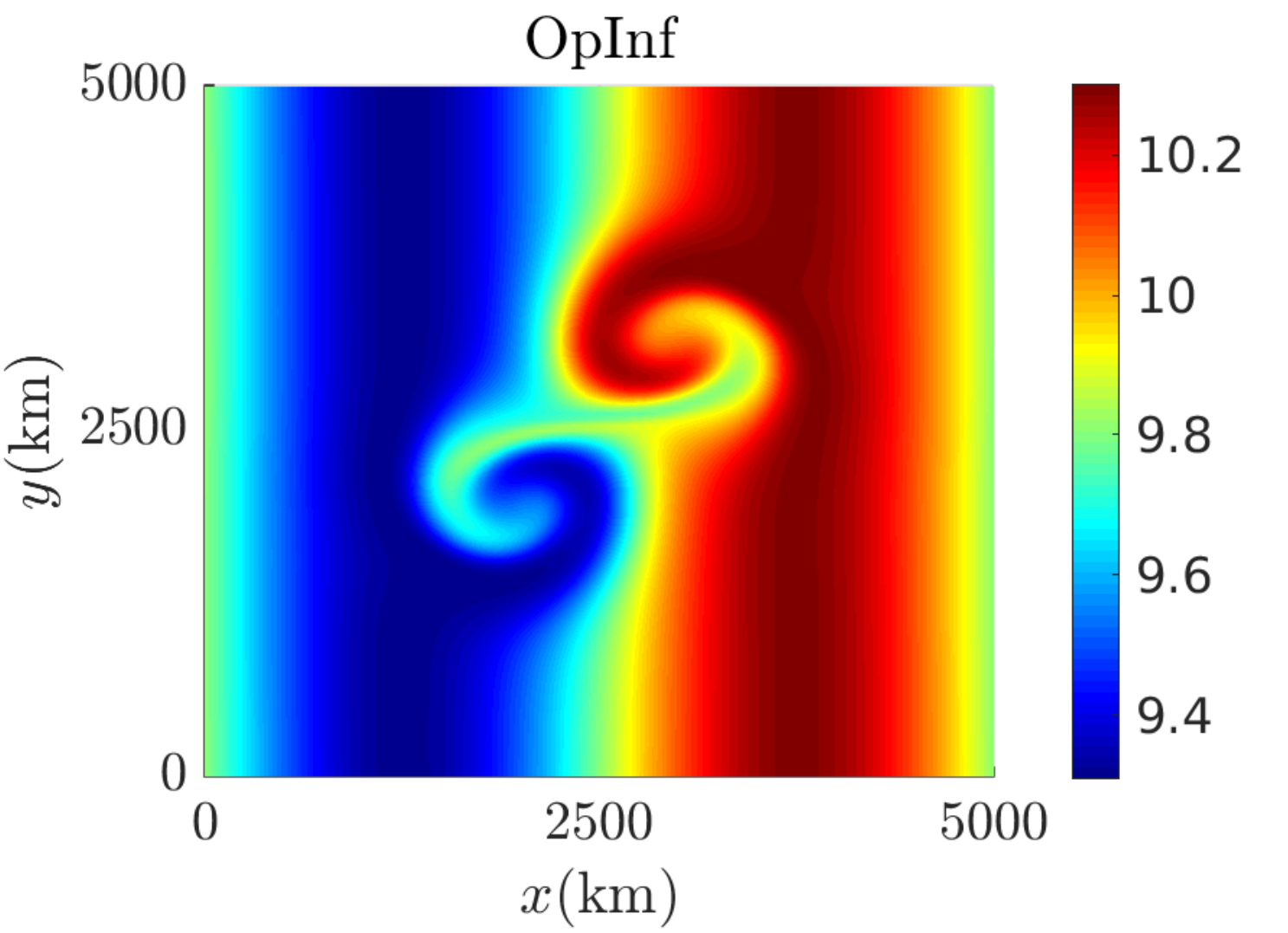}}
	\caption{(Top) Vorticity and (bottom) bouyancy of the FOM and ROMs at the final time. \label{fig:buopar}}	
\end{figure}

In the testing phase we take $ M_{\mathrm test} = 7 $ randomly distributed  parameters
$\mu^{\mathrm test}_{ 1} \ldots, \mu^{\mathrm test}_{M_{\mathrm test}}\in \mathcal{ D}$
that are different from the parameters  in the training phase.
In Figure~\ref{fig:par_rel_error}, we demonstrate the relative errors \eqref{RelAvgError} and \eqref{RelAvgErrortest} for the training and testing periods, respectively, over the increasing number of reduced dimension $r$, which shows that both intrusive and non-intrusive ROMs behave similarly, resulting more accurate reduced solutions with the increasing number of modes. Figure~\ref{fig:buopar} show that the potential vorticity and buoyancy computed at final time for reduced dimension $ r=10 $ with the Coriolis parameter with $\mu^{test}_{ 1}$, are resolved accurately with the intrusive and non-intrusive ROMs.

\section{Conclusions}
\label{sec:conc}

In this paper, we have compared the intrusive POG-G with the data-driven non-intrusive OpInf for a complex fluid dynamics model, i.e., RTSWE. Spatial discretization of the RTSWE results  in a linear-quadratic system in the Hamiltonian form, which is integrated in time with the  energy preserving linearly implicit Kahan's method. In this way, the Hamiltonian structure and the associated invariants of the RTSWE are maintained by the FOM. Numerical results in the parametric and non-parametric settings show that both ROMs behave similarly and yield ROMs that are equal up to numerical errors. This validates the convergence of the learned operators to the intrusively obtained reduced operators under certain conditions on the time discretization in \cite{Peherstorfer16,Peherstorfer20}. The least-squares problem is regularized with the minimum norm solution which was not used in the operator inference context.
Both reduced models are able to accurately re-predict the training data and capture much of the overall system behavior in the prediction period. Due to re-projection, the  OpInf is more costly than the POD-G, nevertheless, speed-up factors of order two are achieved by both ROMs. Moreover, preservation of the conserved quantities by both ROMs indicates the stability of the reduced solutions in long-time integration, which is important for Hamiltonian PDEs like the RTSWE. It also indicates the importance of respecting the physics of the complex problems in ROM application.

Both intrusive and non-intrusive techniques can be compared with machine learning techniques. As a future work, we plan to investigate the reduced models of the SWEs with the Hamiltonian neural networks \cite{Jin20}.

\noindent{\bf Acknowledgemets\/}
This work was supported by 100/2000 Ph.D. Scholarship Program of the Turkish Higher Education Council.


\end{document}